%% file: main.tex
\documentclass[journal,twoside,web]{ieeeconf}
\usepackage{cite}
\usepackage{amsmath,amssymb,amsfonts}
\usepackage{algorithm,algorithmic}
\usepackage{textcomp}
\usepackage{fancyhdr} 
\usepackage{mathtools}
\usepackage{subcaption}

\usepackage{tikz}
\usetikzlibrary{calc}
\usetikzlibrary{positioning}
\usetikzlibrary{arrows.meta}
\usepackage{graphicx}
\usepackage[font={small}]{caption}
\usepackage{setspace}
\usepackage{url}
\usepackage{here}
\usepackage{booktabs}
\usepackage{multirow}
\usepackage{multicol}
\usepackage{bigstrut}
\usepackage{cases}
\usepackage{comment}
\usepackage{framed}
\newtheorem{theorem}{Theorem}
\newtheorem{coro}{Corollary}
\newtheorem{lemma}{Lemma}
\newtheorem{exam}{Example}
\newtheorem{prop}{Proposition}
\newtheorem{rem}{Remark}

\newtheorem{assu}{Assumption}

\newcommand{\CC}{\mathcal{C}}
\newcommand{\CN}{\mathcal{N}}

\newcommand{\CE}{\mathcal{E}}

\newcommand{\CP}{\mathcal{P}}

\newcommand{\BD}{\mathbb{D}}

\newcommand{\QG}{ { \mathcal{Q}_\mathcal{G}^{\mathrm{all}} } }
\newcommand{\QGactive}{\mathcal{Q}_\mathcal{G}}
\newcommand{\QiGactive}{\mathcal{Q}_\mathcal{G}^{i}}
\newcommand{\QjGactive}{\mathcal{Q}_\mathcal{G}^{j}}
\newcommand{\maxQG}{\mathcal{Q}_\mathcal{G}^\mathrm{max}}

\newcommand{\CS}{\mathcal{S}}
\newcommand{\CK}{\mathcal{K}}

\newcommand{\CG}{\mathcal{G}}

\newcommand{\SE}{{E}}

\newcommand{\BR}{\mathbb{R}}

\newcommand{\tlA}{\Tilde{A}}
\newcommand{\tlB}{\Tilde{B}}

\newcommand{\tlP}{\Tilde{P}}
\newcommand{\tlQ}{\Tilde{Q}}
\newcommand{\tlZ}{\Tilde{Z}}
\newcommand{\tlK}{\Tilde{K}}

\newcommand{\tlC}{\Tilde{C}}
\newcommand{\tlD}{\Tilde{D}}
\newcommand{\diag}{\text{diag}}

\newcommand{\bdiag}{\text{blkdiag}}
\newcommand{\He}{\mathrm{He}}

\begin{document}

\title{
Convex Reformulation of LMI-Based Distributed Controller Design \\with a Class of Non-Block-Diagonal Lyapunov Functions
}
\author{{Yuto Watanabe, Sotaro Fushimi, and Kazunori Sakurama}
\thanks{Yuto Watanabe is with the Department of Electrical and Computer Engineering, University of California San
Diego, 9500 Gilman Dr, La Jolla, CA 92093, USA. (email: \texttt{y1watanabe@ucsd.edu})}
\thanks{Sotaro Fushimi is the Department of Mechanical and Aerospace Engineering, University of California San
Diego, 9500 Gilman Dr, La Jolla, CA 92093, USA. (email: \texttt{sfushimi@ucsd.edu})}
\thanks{Kazunori Sakurama is with the Department of System Innovation, Graduate School of Engineering Science, Osaka University, 1-3, Machikaneyama, Toyonaka, Osaka 560-8531, Japan. (email: \tt{sakurama.kazunori.es@osaka-u.ac.jp})}
\thanks{This work was partially supported by the JSPS KAKENHI (under grant 23K22781) and
the joint project of Kyoto University and Toyota Motor Corporation, titled ``Advanced Mathematical Science for Mobility Society''.}
}
\maketitle

\begin{abstract}
This study addresses a distributed state feedback controller design problem for continuous-time linear time-invariant systems by means of linear matrix inequalities (LMIs). As structural constraints on a control gain result in non-convexity in general, the block-diagonal relaxation of Lyapunov functions has been prevalent despite its conservatism. In this work, we target a class of non-block-diagonal Lyapunov functions with the same sparsity pattern as distributed controllers. By leveraging a block-diagonal factorization of sparse matrices and Finsler's lemma, we first present a nonlinear matrix inequality for stabilizing distributed controllers with such Lyapunov functions, which boils down to a necessary and sufficient condition for such controllers if the sparsity pattern is chordal. As its relaxation, we derive novel LMIs, one of which strictly covers the conventional relaxation, and then provide analogous results for $H_\infty$ control. Lastly, numerical examples underscore the efficacy of our results.
\end{abstract}

\section{Introduction}

For the past decades, distributed control has been a pivotal topic in control as a response to the significant advancement of information and communication technologies \cite{siljak2011decentralized,ren2008distributed,bullo2009distributed,molzahn2017survey}.
Its practical applications
include 
power networks \cite{molzahn2017survey,anderson2019system},
formation control \cite{ren2008distributed,bullo2009distributed,yuan2023structured},
process control \cite{schuler2012decentralized}, and traffic systems \cite{halder2020distributed}.
By structuring controllers implemented only with local information, distributed controllers can substantially enhance the scalability and expandability of systems.

Nevertheless, such controller structure makes their design much more challenging even for linear systems \cite{anderson2019system}.
While linear matrix inequalities (LMIs) \cite{boyd1994linear} enable efficient controller design in a convex manner for centralized controllers, 
it is a different story in distributed ones;
the structure constraint for the distributedness results in non-convexity.
As a way to circumvent this issue,
the most prevalent convex restriction is the block-diagonal relaxation of Lyapunov functions \cite{sootla2019existence} as $\sum_{i}x_i^\top P_i x_i$.
However, this relaxation leads to conservatism in general despite the simplicity.
While several works have been proposed to resolve the issue of conservatism 
\cite{apkarian2001continuous,ebihara2004new,pipeleers2009extended,ferrante2019design,furieri2020sparsity,yuan2023structured}, whether an equivalent convex reformulation of the structure constraint for gain matrices exists is still an
open research question, except for specific classes, such as positive systems \cite{tanaka2011bounded,rantzer2015scalable}.
Although the extended LMI approach \cite{apkarian2001continuous,ebihara2004new,pipeleers2009extended} 
gives a dense Lyapunov matrix through the introduction of a slack variable, this approach still requires a block-diagonal relaxation of the slack variable.

For the LMI-based distributed controller design problem,
this study focuses on continuous-time systems and, importantly, a class of non-block-diagonal Lyapunov functions that has the same sparsity pattern as distributed controllers.
For this class of Lyapunov functions,
we propose a new convex relaxation not only for stabilization but also for $H_\infty$ control \cite{scherer2000linear,dullerud2013course}.
Note that such Lyapunov functions play a fundamental role in gradient-flow systems \cite{sakurama2014distributed,sakurama2022generalized}, and
this approach can be easily extended to other fundamental problems, such as $H_2$ control.
First, by leveraging
a block-diagonal factorization for sparse matrices and Finsler's lemma \cite{de2007stability},
we show a new nonlinear matrix inequality for distributed controllers with such Lyapunov functions, which becomes necessary and sufficient over chordal sparsity \cite{vandenberghe2015chordal,zheng2021chordal}.
We then present new LMIs by relaxing the matrix inequality and demonstrate that the $H_\infty$ control versions can be derived similarly.
We also discuss our results in light of a state transformation called \textit{inclusion principle} \cite{siljak2011decentralized,ikeda1982inclusion}, providing an intuitive interpretation.
Finally, numerical examples illustrate the effectiveness of the proposed method.
In the examples, we also show the results of another LMI obtained by slightly modifying the proposed LMI, which has no theoretical guarantee of stabilization but interestingly exhibits excellent performance.

Our contributions can be summarized as follows:
    a) This work presents novel LMIs for the distributed controller design problem.
    They can generate Lyapunov functions with a more complex sparsity pattern than block-diagonal matrices.
    One of the proposed methods completely covers the block-diagonal relaxation. In our numerical experiments, our proposed methods outperform not only the block-diagonal relaxation and sparsity invariance (SI) approaches \cite{furieri2020sparsity} but also extended LMI \cite{apkarian2001continuous,ebihara2004new,pipeleers2009extended};
    b)
    The computation of inverse matrices in
    the proposed method can be decomposed into smaller ones corresponding to subsystems formed by cliques of the communication graph, which alleviates computational burdens and enhances the scalability.
    Note that \cite{ferrante2019design,furieri2020sparsity} do not permit such a decomposed computation of inverse matrices;
    c)
    We show that the derived conditions provide a necessary and sufficient condition for distributed ($H_\infty$) controllers with the class of non-block-diagonal Lyapunov functions under chordal sparsity.
    Since our convex relaxations follow from the condition, one can easily evaluate the conservatism.
    While the authors of \cite{yuan2023structured} utilized a similar block-diagonal factorization, such a discussion cannot be found in any prior work, to our knowledge.
    
    Note that
    our proposed approach can be effectively combined with many existing works, including the extended LMI approach \cite{apkarian2001continuous,ebihara2004new,pipeleers2009extended}.
    By inheriting the advantages of the two methods,
    our combined method enjoys both preferable theoretical and numerical performance.
    For the details,
    see Appendix \ref{Appendix:combination} 
    and the discrete-time counterpart \cite{fushimi2024distributed}.

The remainder of the paper is organized as follows.
Section \ref{sec:problem_statement} provides the target system and problem.
Section \ref{sec:preliminaries} prepares several notions and lemmas for the main results.
In Section \ref{sec:main_result}, we present our main results as a solution to the formulated problem, extending it to $H_\infty$ control and presenting a control-theoretic interpretation.
Section \ref{sec:examples} showcases numerical examples.
Finally, Section \ref{sec:conclusion} concludes this paper.

\paragraph*{Notations}
Let $|\cdot|$ be the number of elements in a countable finite set.
Let $I_n\in \mathbb{R}^{n\times n}$ denote the $n\times n$ identity matrix.
Let $O_{n_1\times n_2}$ be the $n_1\times n_2$ zero matrix.
We omit the subscript when it is obvious.
Let $\mathrm{Im}(E)$ be the image space of the matrix $E$.
Let $\diag(a)$ with $a=[a_1,\ldots,a_n]^\top$ denote the diagonal matrix whose $i$th diagonal entry is $a_i\in\BR$.
Similarly, $\bdiag(\ldots,R_i,\ldots)$ represents the block diagonal matrix consisting of $R_i$.
For $M\in\BR^{m\times n}$ with $\mathrm{rank}(M)<m$, $M^\perp$ represents 
a basis of the null-space if $M$, i.e., $M^\perp$ satisfies
$\{x:Mx=0\} = \mathrm{Im}(M^\perp)$.
For a real square matrix $A\in\BR^{n\times n}$.
the function $\He(\cdot)$ represents $\He(A)=A+A^\top$.

\section{Problem statement}\label{sec:problem_statement}
\subsection{Target system}
Consider a large-scale system consisting of $N$ subsystems.
Suppose that their communication network is modeled by 
a time-invariant undirected graph $\CG = (\CN,\CE)$ with $\CN=\{1,\ldots,N\}$ and an edge set $\CE$.
Namely, the set $\CE$ consists of pairs $(i,j)$ of different nodes $i,j\in\CN$, and we assume $(i,j)\in\CE\Leftrightarrow (j,i)\in\CE$.

Now, consider the following linear time-invariant system:
 \begin{align}\label{eq:LTI_system}
        \dot{x}(t) = 
        Ax(t) + Bu(t) 
    \end{align}
where $x = [x_1^\top,\ldots,x_N^\top]^\top \in \BR^n$ with $x_i\in\BR^{n_i},\,\sum_{i=1}^N n_i = n$, $u = [u_1^\top,\ldots,u_N^\top]^\top \in \BR^m$ with $u_i\in\BR^{m_i},\,\sum_{i=1}^N m_i=m$, $A\in \BR^{n\times n}$, and $B\in\BR^{n\times m}$.
Here, $x_i$ and $u_i$ represent the state and control input of agent $i$.
Besides, $A$ and $B$ are supposed to be partitioned with respect to the partitions $\{n_1,\ldots,n_N\}$ of $n=\sum_{i=1}^N n_i$ and $\{m_1,\ldots,m_N\}$ of $m=\sum_{i=1}^N m_i$ corresponding to the dimension of each subsystem's state and input.

In this paper, we design a static state feedback controller to stabilize the system:
\begin{equation}\label{state_feedback}
    u(t) = K x(t).
\end{equation}
In the following, $K_{ij} \in \BR^{m_i\times n_j}$ represents the $(i,j)$ block of $K$ associated with the partitions $\{n_1,\ldots,n_N\}$ and $\{m_1,\ldots,m_N\}$.

To achieve the stabilization in a distributed fashion,
we impose the following assumption on $K$.
\begin{assu}\label{assumption:K_sparsity}
$K \in \CS$, where
    \begin{equation}\label{def:S}
    \CS = \{K \in \BR^{m\times n}: K_{ij} = O_{m_i\times n_j} \text{ if } (i,j)\notin \CE,\,i\neq j\}.
    \end{equation}
\end{assu}

For $K\in\CS$, we obtain $u_i = \sum_{j=1}^N K_{ij}x_j = \sum_{j\in\CN_i} K_{ij}x_j$ with $\CN_i=\{j\in\CN:(i,j)\in\CE\}\cup\{i\}$, which can be implemented in a distributed fashion.

For simplicity of notations in the following, we assume
$n_i = m_i$ for all $i\in\CN$ without loss of generality;
otherwise, we just need to add the zero vectors to (or remove columns from) $B_i$ to make it square, which is valid because $B_i$ is not necessary to be column full rank in all the results below.

\subsection{Target problem}

For this system, the problem of finding a stabilizing controller satisfying Assumption \ref{assumption:K_sparsity} is equivalent to 
\begin{align}\label{eq:original_problem}
    \text{Find }& P \succ O,\,K\in\CS \nonumber\\
    \text{ s.t. }&
        (A+BK)^\top P + P(A+BK) \prec O.
\end{align}
By implementing \eqref{state_feedback} with $K$ in \eqref{eq:original_problem}, one can guarantee the stability with the Lyapunov function $x^\top P x$.

However, the inequality in \eqref{eq:original_problem} involves non-convexity, which complicates distributed controller design.
Indeed, the standard change of variables, i.e., $Q=P^{-1}$ and $Z = KQ$, gives the following equivalent problem to \eqref{eq:original_problem}:
\begin{align}\label{def:LMI_nonconv}
    \text{Find }& Q \succ O,\,Z \nonumber \\
    \text{ s.t. }&
        QA^\top+AQ + Z^\top B^\top + BZ \prec O,\,
         ZQ^{-1} \in \CS.
\end{align}
Importantly, $K=ZQ^{-1} \in \CS$ is a non-convex constraint,
and the exact convexification is still an open question.

A popular convex relaxation is the following LMI that restricts $Q$ to a block-diagonal matrix, i.e., assumes $x^\top Q^{-1} x = \sum_{i=1}^N x_i^\top Q_i^{-1}x$ as a Lyapunov function, which guarantees $ZQ^{-1}\in\CS \Leftrightarrow Z\in\CS$.
Then, \eqref{def:LMI_nonconv} is reduced to the following relaxed LMI:
\begin{align}\label{LMI_diag}
    \text{Find }& Q=\bdiag(Q_1,\ldots,Q_N) \succ O,\,Z\in\CS \nonumber \\
    \text{ s.t. }&
        QA^\top+AQ + Z^\top B^\top + BZ \prec O,
\end{align}
where $Q_i\in\BR^{n_i\times n_i}$.
Let $\CK_{\CS,\diag}$ be the set of all $K$ given by \eqref{LMI_diag} as follows:
\begin{align*}
\CK_{\CS,\diag} = & \{K=ZQ^{-1}: \exists Q = \bdiag(Q_1,\ldots,Q_N)\succ O, \\&\quad Z\in\CS
    \text{ s.t. } QA^\top +AQ + Z^\top B^\top + BZ \prec O\}.
\end{align*}
Since the matrix $Q$ is non-block-diagonal in general,
this relaxation causes conservatism.

To mitigate the conservatism,
our goal is to find a new subclass of $K$ satisfying \eqref{eq:original_problem} that generalizes $\CK_{\CS,\diag}$ and can be obtained in a convex manner, e.g., via LMIs.
For this purpose, we here focus on the following problem and give a solution in the form of LMI.
\paragraph*{Problem 1}
Consider linear time-invariant system \eqref{eq:LTI_system} with a static state feedback controller in \eqref{state_feedback} over an undirected graph $\CG$ and the set $\CS$ in \eqref{def:S}. Then, solve the following:
\begin{align}\label{problem1}
        \text{Find} & \quad
        P\succ O,\,P\in\CS,\, K\in\CS \nonumber\\
        \text{ s.t. }&
        \quad 
        (A+BK)^\top P + P(A+BK) \prec O.
\end{align}

In the following, we let $\CK_\CS$ be the set of all $K$ satisfying \eqref{problem1}, i.e.,
\begin{align}\label{def:LMI}
    \CK_\CS =& \{K\in\CS:
    \exists P\in\CS,\,P\succ O \text{ s.t. } \eqref{problem1}
    \}.
\end{align}
Although \eqref{problem1} still involves a relaxation owing to $P\in\CS$ compared with the original problem \eqref{def:LMI_nonconv},
this can generalize the conventional block-diagonal relaxation in \eqref{LMI_diag}, i.e., $\CK_{\CS,\diag}\subset\CK_\CS$.
Note that for gradient-flow systems,
this class of Lyapunov functions characterizes all quadratic potential functions that generate a distributed gradient \cite{sakurama2014distributed}.
In addition, even if $n_i\neq m_i$, we can consider the same sparsity pattern for $K$ by setting an appropriate dimension.

To take the robustness and exogenous disturbances into account,
we also tackle $H_\infty$ control in Subsection \ref{subsec:H_2_infty}.

\section{Preliminaries}\label{sec:preliminaries}

Before establishing our main results in Section \ref{sec:main_result},
let us introduce several notions and supporting lemmas.

\paragraph{Graph theory}
First, we prepare graph theoretic concepts.
Consider an undirected graph $\CG=(\CN,\CE)$.
A \textit{clique} $\CC$ is a subset of $\CN$ satisfying $(i,j)\in\CE$ for all $i,j\in\CC,\,i\neq j$.
Namely, a clique is a node set of a complete subgraph of $\CG$.
We assign indices to the cliques of $\CG$ as $\CC_1,\CC_2,\ldots$ and let $\QG$ denote the set of all cliques' indices.
If $\CC$ is not included in any other cliques, $\CC$ is said to be \textit{maximal}.
Let $\maxQG (\subset \QG)$ represent the index set of all maximal cliques of $\CG$.
In the following, we use $\QGactive$ as a subset of $\QG$ (not necessarily as $\maxQG$) and $\QiGactive$ as the subset of cliques in $\QGactive$ by which node $i$ is contained, i.e., $\QiGactive:=\{k\in\QGactive: i\in\CC_k\}$.
A graph $\CG$ is said to be \textit{chordal} if every cycle of length at least four has a chord \cite{vandenberghe2015chordal,zheng2021chordal}.
Note that this class is not restrictive, and
readers can find examples of chordal graphs in \cite{vandenberghe2015chordal,zheng2021chordal}.

\paragraph{Matrix $\SE$}
Consider a subset $\QGactive$ of $\QG$.
For a clique $\CC_k,\,k\in\QGactive$, we define $\SE_{\CC_k}\in \BR^{n_{\CC_k}\times n}$ with $n_{\CC_k}=\sum_{j\in\CC_k}n_j$ as
\begin{equation}\label{E_Ck}
    \SE_{\CC_k} = 
    \underset{j\in \CC_k}
    {
    [\ldots,\SE_j^\top,\ldots]^\top,
    }
\end{equation}
where $\SE_j = [O_{n_j\times n_1},\ldots,I_{n_j},\ldots,O_{n_j\times n_N}] \in \BR^{n_j\times n}$.
This matrix generates the clique-wise copy of $x=[x_1^\top,\ldots,x_N^\top]^\top$ as 
    $\SE_{\CC_k}x = \underset{j\in \CC_k}
    {[\ldots,x_j^\top,\ldots]^\top}$.

Interestingly, for a subset $\QGactive$ of $\QG$,
the matrices $\SE_{\CC_k},\,k\in\QGactive$ have the following beneficial properties for distributed controller design.
The proof can be found in \cite[Lemma 2]{watanabe2023distributed}.
\begin{prop}\label{prop:CD_matrix}
Consider the following matrix $\SE$ consisting of $\SE_{\CC_k}$ in \eqref{E_Ck}:
\begin{equation*}
    \SE = \underset{k\in\QGactive}{[\ldots,\SE_{\CC_k}^\top,\ldots]^\top} \in \BR^{(\sum_{k\in\QGactive}n_{\CC_k})\times n}.
\end{equation*}
Assume that 
$\QiGactive\neq\emptyset$ for all $i\in\CN$.
Then, the matrix $\SE$ satisfies the following properties.
\begin{enumerate}
    \item[(a)] $\SE$ is column full rank.
    \item[(b)] $\SE^\top \SE = \bdiag(|\QGactive^1|I_{n_1},\ldots,|\QGactive^N|I_{n_N}) \succ O$.
    \item[(c)] For $\Tilde{x} = [\ldots,\Tilde{x}_k^\top,\ldots]^\top \in \BR^{\sum_{k\in\QGactive}n_{\CC_k}}$ with $\Tilde{x}_k\in\BR^{n_{\CC_k}}$,
    \begin{equation}\label{eq:E_transpose}
        \SE^\top \Tilde{x} = 
        \begin{bmatrix}
            \sum_{k\in\QGactive^1} E_{k,1} \Tilde{x}_k\\
            \vdots\\
            \sum_{k\in\QGactive^N} E_{k,N} \Tilde{x}_k
        \end{bmatrix} \in \BR^{n},
    \end{equation}
    where $E_{k,i}\in\BR^{n_i\times n_{\CC_k}}$ for $k\in\QGactive$ and $i\in\CC_k$ satisfies $E_{k,i} x_{\CC_k}=x_i$ with $x_{\CC_k} = \underset{j\in\CC_k}{[\ldots,x_{j}^\top,\ldots]^\top} \in \BR^{n_{\CC_k}}$ for any $x=[x_1^\top,\ldots,x_N^\top]^\top\in\BR^n$.
\end{enumerate}
\end{prop}

This proposition together with $\bigcup_{k\in\QiGactive}\CC_k\subset\CN_i$ \cite{sakurama2022generalized} implies that one can compute the least squares solution $(\SE^\top\SE)^{-1}\SE^\top \Tilde{x}$ of $\Tilde{x}=\SE x$ and
the projection $\SE(\SE^\top\SE)^{-1}\SE^\top\Tilde{x}$ onto $\mathrm{Im}(\SE)$ in a distributed fashion.
Additionally, the assumption $\QiGactive\neq \emptyset,\,i\in\CN$, which means that each node belongs to at least one clique, is always satisfied for $\QGactive=\maxQG,\,\QG$, regardless of the connectivity of $\CG$.

\begin{exam}
Consider a system with $N=3$ over the chordal graph in Fig. \ref{fig:example} with maximal cliques $\QGactive=\maxQG=\{1,2\}$ with $\CC_1=\{1,2\}$ and $\CC_2=\{2,3\}$.
Then, $\QGactive^1=1$,  $\QGactive^2=2$, and  $\QGactive^3=1$.
When $n_1=n_2=n_3 =1$, the $\SE$ matrix is given by $\SE = [\SE_{\CC_1}^\top,\,\SE_{\CC_2}^\top]^\top \in\BR^{4\times 3}$ with
\begin{equation*}
    \SE_{\CC_1} = 
    \footnotesize{\begin{bmatrix}
    1 & 0 & 0\\
    0 & 1 & 0
    \end{bmatrix}}, \quad
     \SE_{\CC_2} = 
    \footnotesize{\begin{bmatrix}
    0 & 1 & 0\\
    0 & 0 & 1
    \end{bmatrix}}.
\end{equation*}
This $\SE$ matrix satisfies $\SE_{\CC_1} x = [x_1,x_2]^\top$,  $\SE_{\CC_2} x = [x_2,x_3]^\top$ for $x=[x_1,x_2,x_3]^\top\in\BR^3$, and $\SE^\top\SE = \diag(1,2,1) = \diag(|\QGactive^1|,|\QGactive^2|,|\QGactive^3|)$.
\begin{figure}[t]
    \centering
    \includegraphics[width=0.75\columnwidth]{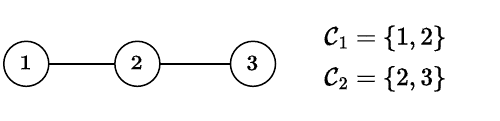}
    \caption{An example of systems with $N=3$. This network is chordal with maximal cliques $\CC_1=\{1,2\}$ and $\CC_2=\{2,3\}$.}
    \label{fig:example}
\end{figure}
\end{exam}
\begin{rem}
The matrix $\SE$ is called the \textit{clique-wise duplication (CD) matrix} in our previous work \cite{watanabe2023distributed} in the context of distributed optimization, which provides more detailed properties of this matrix.
Moreover, there is more flexibility in the dimension of each block in $E$; see \cite[Appendix A]{fushimi2024distributed}.
\end{rem}
\paragraph{Block-diagonal factorization of sparse matrices}
Leveraging the matrix $\SE$, we first present the positive definite version of Agler's theorem \cite{grone1984positive,vandenberghe2015chordal,zheng2021chordal}.
This theorem ensures the existence of a block-diagonal factorization of sparse positive definite matrices over chordal graphs.
\begin{lemma}[{positive definite version of Agler's theorem}]\label{lemma:agler}
    Consider undirected graph $\CG$ with cliques $\QGactive=\{1,\ldots,q\}$.
    Then, given a partition $\{n_1,\ldots,n_N\}$ of $n$,
    for $\tlP = \bdiag(\tlP_1,\ldots,\tlP_q)\succ O
    $ with $\tlP_k\in\BR^{n_{\CC_k}\times n_{\CC_k}}$, the matrix $P = \SE^\top \tlP\SE = \sum_{k=1}^q \SE_{\CC_k}^\top \tlP_k \SE_{\CC_k}$ is positive definite and belongs to $\CS$.
    Moreover,
    if $\CG=(\CN,\CE)$ is chordal with maximal cliques $\QGactive=\maxQG=\{1,\ldots,q\}$, the following equivalence holds:
        \begin{align}\label{eq:agler}  
    &P \in \CS \text{ and } P\succ O \nonumber \\
    \Leftrightarrow&\exists 
        \tlP = \bdiag(\tlP_1,\ldots,\tlP_q)\succ O
        \text{ with }\tlP_k\in\BR^{n_{\CC_k}\times n_{\CC_k}}\nonumber \\
        &\text{ s.t. }
        P = \SE^\top \tlP\SE = \sum_{k=1}^q \SE_{\CC_k}^\top \tlP_k \SE_{\CC_k}.
    \end{align}
\end{lemma}

Next, for all undirected graphs (not necessarily chordal),
we show that every matrix $K$ in $\CS$ has the following block-diagonal factorization. This lemma implies that $K\in\CS$ can be parameterized by a clique-wise block-diagonal matrix $\tlK$.
\begin{lemma}\label{lemma:sparsity}
    Consider a partition $\{n_1,\ldots,n_N\}$ of $n$.
    Let $\QiGactive,\,i\in\CN$ satisfy $\QiGactive\neq \emptyset$. Suppose $\QiGactive\cap \QjGactive \neq \emptyset \Leftrightarrow (i,j)\in\CE$.
    Then, it holds that
     \begin{equation*}
            \CS = \{\SE^\top \tlK \SE : \tlK = \bdiag(\underset{k\in\QGactive}{\cdots,\tlK_k,\cdots}),\,\tlK_k\in \BR^{n_{\CC_k}\times n_{\CC_k}}
            \}.
        \end{equation*}
\end{lemma}
\begin{proof}
See Appendix \ref{Appendix:proof-sparsity}.
\end{proof}

\begin{exam}
For the graph in Fig. \ref{fig:example}, Lemma \ref{lemma:sparsity} can be verified as follows.
Here $*_{2\times 2}$ is an arbitrary $2\times 2$ matrix.
\begin{equation*}
\SE^\top 
\begin{bmatrix}
    *_{2\times 2}& O_{2\times 2}\\
    O_{2\times 2}& *_{2\times 2}\\
\end{bmatrix}
\SE =
\footnotesize
{\begin{bmatrix}
* & * & 0\\
* & * & *\\
0 & * & *
\end{bmatrix}} \in \CS.
\end{equation*}
\end{exam}

\section{Main results}\label{sec:main_result}

This section presents our main results for Problem 1.
The notion of chordal graphs and Lemmas \ref{lemma:agler} and \ref{lemma:sparsity}
prepared in Section \ref{sec:preliminaries} allow us to derive a (nonlinear) matrix inequality to solve Problem 1, which yields novel LMIs.
Over chordal graphs, this nonlinear matrix inequality boils down to a necessary and sufficient condition for $K\in\CK_\CS$.
Additionally, we apply the same strategy to the $H_\infty$ control problem, presenting analogous LMIs.
Note that our proposed approach can be extended to other important scenarios, such as polytopic uncertainties and output feedback control.
For the details,
see Appendix \ref{Appendix:extention}.

In what follows, let $\QGactive=\{1,\ldots,q\}$ without loss of generality, and
for notational simplicity, define the following sets of block-diagonal matrices corresponding to $\QGactive$:
\begin{align*}
    &\BD_{\QGactive} = \{
    \bdiag(\tlP_1,\ldots,\tlP_q):
    \tlP_k \in \BR^{n_{\CC_k}\times n_{\CC_k}},\,\forall k\in\QGactive
    \},\\
    &\BD_{\QGactive}^{++} = \{
    \bdiag(\tlP_1,\ldots,\tlP_q)
    \in \BD_{\QGactive}:
    \tlP_k\succ 0,\,\forall k\in\QGactive
    \}.
\end{align*}
\subsection{Solution to Problem 1}\label{subsec:solution_to_problem_1}

A solution to Problem 1 is presented in this subsection.
First, we impose the following assumption which is not strict as it is always satisfied for $\QGactive=\QG$ and $\maxQG$, regardless of the connectivity of $\CG$.
\begin{assu}\label{assumption:Q}
For $\CG$,
the index set $\QGactive=\{1,\ldots,q\}$ of its cliques satisfies the following:
\begin{itemize}
    \item $\QiGactive\neq \emptyset$ for all $i\in\CN$.
    \item $\QiGactive\cap \QjGactive \neq \emptyset \Leftrightarrow (i,j)\in\CE$.
\end{itemize}
\end{assu}

Next, the following lemma, known as \textit{Finsler's lemma}, plays a key role together with the lemmas in Section \ref{sec:preliminaries}.
\begin{lemma}[Finsler's lemma \cite{de2007stability}]
\label{lemma:finsler}
    Let $x\in\BR^n$, $Q=Q^\top \in \BR^{n\times n}$, and $M\in\BR^{r\times n}$ such that $\mathrm{rank}(M)<n$.
    The following statements are equivalent:
    \begin{enumerate}
    \item $x^\top Q x <0,\, \forall Mx=0,\,x\neq 0$.
    \item $M^{\perp\top} Q M^{\perp} \prec O$.
    \item $\exists\rho\in\BR$ s.t. $Q + \rho M^\top M\prec O$.
    \item $\exists X\in \mathbb{R}^{n\times r}$ s.t. $Q+M^\top X^\top + XM \prec O$.
    \end{enumerate}
\end{lemma}

 In our main result,
 we use $I-\SE(\SE^\top\SE)^{-1}\SE^\top$ 
 as the matrix $M$ in Finsler's lemma based on the following proposition, which follows from Proposition \ref{prop:CD_matrix}.
\begin{prop}\label{prop:perp}
Consider undirected graph $\CG=(\CN,\CE)$ with a set of cliques $\QGactive\subset\QG$.
Assume that $\CG$ is not complete.
Suppose $\QiGactive\neq\emptyset$ for all $i\in\CN$.
Consider 
\begin{equation}\label{eq:M}
    M=I-\SE(\SE^\top\SE)^{-1}\SE^\top.
\end{equation}
Then, $\mathrm{rank}(M)<\sum_{k\in\QGactive} n_{\CC_k}$ and $M^\perp = \SE$ hold.
\end{prop}
\begin{proof}
See Appendix \ref{Appendix:perp}.
\end{proof}

As a preliminary for controller design, let us present the following lemma for the internal stability, which characterizes $P\in\CS$ in Problem 1 with $K =O$ if $\CG$ is chordal.
For non-chordal graphs, this lemma still gives a sufficient condition for such a $P\in\CS$.
Note that for the system matrices $A,\,B$ in \eqref{eq:LTI_system}, the matrices $\tlA,\,\tlB$ below represent
\begin{align}\label{def:extended_matrices}
  &\tlA = \SE A (\SE^\top\SE)^{-1} \SE^\top,\, 
  \tlB = \SE B (\SE^\top\SE)^{-1} \SE^\top,
\end{align}
whose interpretation is presented in Subsection \ref{subsec:inclusion_principle}.

\begin{lemma}\label{lemma:stability_analysis}
Suppose Assumption \ref{assumption:Q}.
Let 
\begin{align*}
\CP_\CS =& \{P\in\CS:
    P\succ O,\, A^\top P + PA \prec O
    \} \\
\hat{\CP}_\CS =& \{P= \SE^\top \tlP \SE 
:\exists\tlP\in\BD_{\QGactive}^{++},\,
    \rho\in\BR\\
    &
    \text{s.t. } 
    \tlA^\top \tlP + \tlP\tlA +
    \rho M
    \prec O
    \}
\end{align*}
with $M$ in \eqref{eq:M}.
Then, $\hat{\CP}_\CS \subset\CP_\CS$ holds.
Moreover, $ \hat{\CP}_\CS  = \CP_\CS$
holds if $\CG$ is chordal and $\QGactive=\maxQG$.
\end{lemma}
\begin{proof}
When $\CG$ is complete, then $\SE$ is reduced to the identity matrix and thus we have $M=O$. Therefore $\CP_\CS$ and $\hat{\CP}_\CS$ are equivalent.

It is assumed below that $\CG$ is not complete.
First, $\hat{\CP}_\CS \subset \CP_\CS$ is shown as follows.
Consider $\SE^\top \tlP \SE \in\hat{\CP}_\CS$. Then, we obtain
\begin{align*}
&\SE^\top (\tlA^\top \tlP + \tlP\tlA +
    \rho (I - \SE(\SE^\top\SE)^{-1}\SE^\top) ) \SE \\
= & A^\top (\SE^\top \tlP \SE) + (\SE^\top\tlP\SE)A^\top \prec O.
\end{align*}
Since Lemma \ref{lemma:sparsity} guarantees $\SE^\top\tlP\SE\in\CS$, we have $\hat{\CP}_\CS \subset \CP_\CS$.

To prove $\hat{\CP}_\CS=\CP_\CS$ for chordal graph $\CG$ with $\QGactive=\maxQG$, we show $\CP_\CS\subset \hat{\CP}_\CS$.
Let $P\in\CP_\CS$.
Then, by Lemma \ref{lemma:agler},
$P$ always admits the factorization $P = \sum_{k=1}^q \SE_{\CC_k}^\top \tlP_k \SE_{\CC_k} = \SE^\top \tlP \SE$ as $\CG$ is chordal,
where $\tlP\in \BD_{\QGactive}^{++}$.
Then, by leveraging this, we obtain
\begin{align*}
&A^\top P + PA 
= A^\top (\SE^\top \tlP\SE) +( \SE^\top \tlP\SE) A\\
&= \SE^\top( \SE(\SE^\top\SE)^{-1}A^\top \SE^\top \tlP
+ \tlP \SE A (\SE^\top\SE)^{-1}\SE^\top) \SE \\
&= \SE^\top (\tlA^\top \tlP + \tlP\tlA) \SE \prec O.
\end{align*}
By Finsler's lemma,
the above inequality is equivalent to 
    $\tlA^\top \tlP + \tlP\tlA + \rho M
    \prec O$
for some $\rho\in\BR$.
Therefore $P\in\hat{\CP}_\CS$ holds.
\end{proof}

We now proceed to distributed controller design.
As the first main result,
we present the following theorem as a sufficient condition for distributed controllers $K\in\CK_\CS$.
Furthermore, this condition satisfies the necessity over chordal graphs, meaning that this provides a characterization of the solutions to Problem 1.
Recall that $E^\top E$ is diagonal from Proposition \ref{prop:CD_matrix}, which guarantees $(E^\top E)^{-1}F\in\CS$ for $F\in\CS$.

\begin{theorem}\label{theorem:controller_LMI}
Suppose Assumption \ref{assumption:Q}.
Let
\begin{align*}
&\hat{\CK}_\CS = \{ K=(\SE^\top\SE)^{-1}\SE^\top (\tlZ\tlQ^{-1}) \SE:\\
    &\exists \tlZ \in\BD_{\QGactive},\,
    \tlQ\in\BD_{\QGactive}^{++}
    ,\,
    \rho\in\BR
    \text{ s.t. } \Phi(\tlQ,\tlZ)+\rho \tlQ M\tlQ \prec O
    \}
\end{align*}
with $M$ in \eqref{eq:M} and
\begin{equation}\label{eq:Phi}
\Phi(\tlQ,\tlZ)= \tlQ\tlA^\top + \tlA\tlQ + \tlZ^\top\tlB^\top + \tlB\tlZ.
\end{equation}
Then, $\hat{\CK}_\CS \subset \CK_\CS$ holds for $\CK_\CS$ in \eqref{def:LMI}.
Moreover,
    $\hat{\CK}_\CS = \CK_\CS$
holds if $\CG$ is chordal and $\QGactive=\maxQG$.
\end{theorem}
\begin{proof}
    If $\CG$ is complete, both $\hat{\CK}_\CS$ and $\CK_\CS$ coincide with the set of all the stabilizing controllers because $\SE=I$ and $M=O$. Therefore, the statement is true.
    
    In what follows, we assume that $\CG$ is not complete.
    First, we prove $\hat{\CK}_\CS \subset \CK_\CS$.
    For $(\SE^\top\SE)^{-1}\SE^\top \tlK \SE \in \hat{\CK}_\CS$ with $\tlK=\tlZ\tlQ^{-1}$,
    it can be seen that
    \begin{align*}
    O \succ& \tlQ^{-1}(\tlQ\tlA^\top + \tlA\tlQ + \tlZ^\top\tlB^\top + \tlB\tlQ +\rho \tlQ M \tlQ)\tlQ^{-1}\\
    =& (\tlA+\tlB \tlK)^\top \tlQ^{-1} +
    \tlQ^{-1}(\tlA+\tlB \tlK) + \rho M.
    \end{align*}
    Recalling Finsler's lemma and defining $K = (\SE^\top\SE)^{-1}\SE^\top \tlK \SE$ and $P = \SE^\top \tlQ^{-1}\SE$ yield the following equivalent inequality:
    \begin{align*}
        O\succ
        &\SE^\top( (\tlA+\tlB \tlK)^\top \tlQ^{-1} +
            \tlQ^{-1}(\tlA+\tlB \tlK) ) \SE \\
        =&(A+ B K)^\top P +
        P(A+ BK).
    \end{align*}
    Since Lemma \ref{lemma:sparsity} ensures $K\in\CS$ and $P\in\CS$,
    we obtain $\hat{\CK}_\CS \subset \CK_\CS$.

   Next, let us move on to the proof of $\hat{\CK}_\CS = \CK_\CS$ under chordal graphs with $\QGactive=\maxQG$.
    We shall show $\CK_\CS\subset\hat{\CK}_\CS$.
    Consider $K\in\CK_\CS$.
    By Lemma \ref{lemma:sparsity}, 
    $K$ can be represented by $K = (\SE^\top\SE)^{-1}\SE^\top \tlK \SE$ with $\tlK = \bdiag(\tlK_1,\ldots,\tlK_q) \in \BD_{\QGactive}$, where $(\SE^\top\SE)^{-1}$ is diagonal. 
    Moreover, the existence of a positive definite Lyapunov matrix $P\in\CS$ is guaranteed.
    Then, leveraging Lemma \ref{lemma:agler} for $P \in\CS$,
    we obtain $P = \SE^\top \tlP\SE$ with some $\tlP=\bdiag(\tlP_1,\ldots,\tlP_q)\succ O$.
    Consequently, Lemma \ref{lemma:stability_analysis} gives the following inequality with some $\rho\in\BR$:
    \begin{equation*}
        (\tlA + \tlB\tlK)^\top\tlP
         +\tlP(\tlA + \tlB\tlK) + \rho M \prec O.
    \end{equation*}
    Then, letting $\tlQ = \tlP^{-1}$ and $\tlZ = \tlK\tlQ$ gives the equivalent inequality below:
    \begin{equation*}
        \tlQ \tlA^\top + \tlA\tlQ
        + \tlZ^\top\tlB^\top + \tlB\tlZ
        +\rho \tlQ M \tlQ \prec O.
    \end{equation*}
    Hence, we obtain $K=(\SE^\top\SE)^{-1}\SE^\top \tlK \SE \in \hat{\CK}_\CS$.
    Therefore we arrive at $\hat{\CK}_\CS = \CK_\CS$.
\end{proof} 

Theorem \ref{theorem:controller_LMI} gives the following novel convex relaxation that strictly contains the block-diagonal relaxation $\CK_{\CS,\diag}$.
A solution to the LMI in $\hat{\CK}_\CS$ admits $x\SE^\top \tlQ^{-1}\SE x$ as a Lyapunov function, where $\SE^\top \tlQ^{-1}\SE\in\CS$ is satisfied from Lemma \ref{lemma:sparsity}.
\begin{theorem}[Proposed method 1]\label{theorem:LMI_1}
Suppose Assumption \ref{assumption:Q}.
Let 
\begin{align*}
&\hat{\CK}_{\CS,\mathrm{rlx}1} =  \{ (\SE^\top\SE)^{-1}\SE^\top (\tlZ\tlQ^{-1}) \SE
    :
    \exists\tlZ\in\BD_{\QGactive},\,
    \tlQ\in\BD_{\QGactive}^{++},\\
    &\rho\in\BR,\,\eta>0 
    \text{ s.t. } 
    \Phi(\tlQ,\tlZ)
    + \rho M
     \prec O,\,
    \tlQ M + M \tlQ \succeq \eta M\}
\end{align*}
with $\Phi(\tlQ,\tlZ)$ in \eqref{eq:Phi}.
Then, the following inclusion holds:
\begin{equation*}
    \CK_{\CS,\diag} \subset \hat{\CK}_{\CS,\mathrm{rlx}1}
     \subset \hat{\CK}_{\CS}
    \subset \CK_\CS.
\end{equation*}
Moreover, $\CK_{\CS,\diag} \subset \hat{\CK}_{\CS,\mathrm{rlx}1}
     = \hat{\CK}_{\CS}
    = \CK_\CS$ if $\CG$ is complete.
\end{theorem}
\begin{proof}
First, let us show $\hat{\CK}_{\CS,\mathrm{rlx}1} \subset \hat{\CK}_\CS$.
Consider $K=(\SE^\top\SE)^{-1}\SE^\top (\tlZ\tlQ^{-1}) \SE\in\hat{\CK}_{\CS,\mathrm{rlx}1}$.
Lemma \ref{lemma:sparsity} guarantees $K\in\CS$.
Then, combining the two inequalities in $\hat{\CK}_{\CS,\mathrm{rlx}1}$ gives
$\Phi(\tlQ,\tlZ) \prec \frac{|\rho|}{\eta}(\tlQ M + M\tlQ)$.
Thus, pre- and post-multiplying this by $\SE^\top \tlQ^{-1}$ and $\tlQ^{-1}\SE$ respectively and utilizing Finsler's lemma yield
\begin{align*}
&O\succ \SE^\top \tlQ^{-1}\Phi(\tlQ,\tlZ)\tlQ^{-1}\SE\\ 
{\Leftrightarrow} &\; 
\exists\beta\in\BR \text{ s.t. }
\tlQ^{-1}\Phi(\tlQ,\tlZ)\tlQ^{-1} + \beta M \prec O\\
\Leftrightarrow&\;
\exists\beta\in\BR \text{ s.t. }
\Phi(\tlQ,\tlZ) + \beta \tlQ M\tlQ \prec O.
\end{align*}
Therefore, $\hat{\CK}_{\CS,\mathrm{rlx}1} \subset \hat{\CK}_\CS$.

Next, we prove $\CK_{\CS,\diag} \subset \hat{\CK}_{\CS,\mathrm{rlx}1}$.
Consider $K\in\CK_{\CS,\diag}$.
    Then, we have $K=ZQ^{-1}$ with some $Z\in\CS$ and block-diagonal matrix $Q=\{Q_1,\ldots,Q_N\}\succ O$.
    Additionally, Lemma \ref{lemma:sparsity} guarantees that there exists a block-diagonal matrix $\tlK$ such that $K=(\SE^\top\SE)^{-1}\SE^\top\tlK\SE$.
    Now, by setting $\tlQ =\bdiag(\cdots,\tlQ_k,\cdots)\succ O$ with
    \begin{equation*}
        \tlQ_k = \bdiag(\cdots,\underset{j\in\CC_k}{
        |\QjGactive| Q_j}, \cdots)\succ O
    \end{equation*}
    and $\tlZ = \tlK\tlQ$,
    we have $(\SE^\top\SE)^{-1}\SE^\top \tlZ\tlQ^{-1} \SE =(\SE^\top\SE)^{-1}\SE^\top \tlZ\SE(\SE^\top\SE)^{-1} Q^{-1} = ZQ^{-1}=K$ since $\SE (\SE^\top\SE)^{-1}Q^{-1}= \tlQ^{-1}\SE$
    from the proof of Proposition \ref{prop:diag-Q-E} in Appendix \ref{Appendix:perp}.
    Note that $Z$ can be represented as $Z = (\SE^\top\SE)^{-1}\SE^\top \tlZ\SE(\SE^\top\SE)^{-1}\in\CS$ and it holds for $\tlQ$ and $M$ that
\begin{equation}\label{eq:MQ=QM}
    \tlQ M = M\tlQ = \tlQ - \SE^\top Q \SE
\end{equation}
from Proposition \ref{prop:diag-Q-E}.
Hence, we obtain
        \begin{align*}
        O&\succ (A+BK)^\top Q^{-1} + Q^{-1}(A+BK)\\
        &= (A+BK)^\top (\SE^\top \tlQ^{-1}\SE) + (\SE^\top \tlQ^{-1}\SE) (A+BK)\\
        &= \SE^\top ((\tlA + \tlB \tlK)^\top \tlQ^{-1}
        + \tlQ^{-1} (\tlA+\tlB\tlK)) \SE.
    \end{align*}
    Thus, by using Finsler's lemma and following the same procedure as the proof $\CK_\CS\subset\hat{\CK}_\CS$ over chordal graphs above, one can obtain $K\in\hat{\CK}_\CS$, which gives $\CK_{\CS,\diag}\subset\hat{\CK}_\CS$.
Then, utilizing Theorem \ref{theorem:controller_LMI} and \eqref{eq:MQ=QM}, we get $O\succ\Phi(\tlQ,\tlZ)+\nu \tlQ M \tlQ
= \Phi(\tlQ,\tlZ) + \nu M\tlQ^2 M
\succeq  \Phi(\tlQ,\tlZ) + \rho M$ for some $\nu\in\BR$ and $\rho = \nu \lambda^2 \in\BR$, where $\lambda$ represents either the maximal or minimum eigenvalue of $\tlQ$.
Moreover, \eqref{eq:MQ=QM} also gives
\begin{equation*}
M\tlQ + \tlQ M = M^2 \tlQ + \tlQ M^2
= 2 M\tlQ M \succeq \eta M,
\end{equation*}
where $\eta = 2\lambda_{\min}(\tlQ)>0 $ and
$\lambda_{\min}(\tlQ)$ is the minimum eigenvauel of $\tlQ$.
Therefore, we obtain $K\in\hat{\CK}_{\CS,\mathrm{rlx}1}$, and thus $\CK_{\CS,\diag}\subset\hat{\CK}_{\CS,\mathrm{rlx}1}$.
Hence, combining the above discussion with Theorem \ref{theorem:controller_LMI}, we obtain $ \CK_{\CS,\diag} \subset \hat{\CK}_{\CS,\mathrm{rlx}1}
     \subset \hat{\CK}_{\CS}
    \subset \CK_\CS$.
    
If $\CG$ is complete, we have $\SE=I$ and $M=O$, which obviously yields $ \CK_{\CS,\diag} \subset \hat{\CK}_{\CS,\mathrm{rlx}1}
     = \hat{\CK}_{\CS}
    = \CK_\CS$.
\end{proof}

Furthermore, we present another convex restriction that imposes $\rho=0$ on $\hat{\mathcal{K}}_\CS$.
Note that there is no extra constraint as $\tlQ M+\tlQ M \succeq \eta M$ in this variant, which enables searching solutions that cannot be obtained by
the proposed method 1.
\begin{coro}[Proposed method 2]\label{coro:LMI_2}
Suppose Assumption \ref{assumption:Q}.
Let 
\begin{align*}
\hat{\CK}_{\CS,\mathrm{rlx}2} = & \{ (\SE^\top\SE)^{-1}\SE^\top (\tlZ\tlQ^{-1}) \SE
    :
    \\&
    \exists\tlZ\in\BD_{\QGactive},\,
    \tlQ \in\BD_{\QGactive}^{++}
    \text{ s.t. }  
    \Phi(\tlQ,\tlZ)\prec O
    \}
\end{align*}
with $\Phi(\tlQ,\tlZ)$ in \eqref{eq:Phi}.
Then, $\hat{\CK}_{\CS,\mathrm{rlx}2}
    \subset \hat{\CK}_\CS$ holds.
Moreover, $\CK_{\CS,\diag} \subset \hat{\CK}_{\CS,\mathrm{rlx}2}
     = \hat{\CK}_{\CS}
    = \CK_\CS$ if $\CG$ is complete.
\end{coro}
\begin{rem}\label{rem:rho=0}
    In
    the LMI in Corollary \ref{coro:LMI_2},
    we restrict $\rho$ to be zero.
    A similar approach can be found in \cite{yuan2023structured}.
    Note that in general, we cannot expect $\CK_{\CS,\mathrm{diag}}\subset \hat{\CK}_{\CS,\mathrm{rlx}2}$ because this approach cannot cover the case where the optimal value of $\rho$ is negative.
\end{rem}
\begin{rem}[Proposed method 3]\label{rem:LMI_without_stab}
Removing the constraint $\tlQ M + M\tlQ \succeq \eta M,\,\eta>0$ from $\hat{\CK}_{\CS,\mathrm{rlx}1}$ in Theorem \ref{theorem:LMI_1} gives another set of $K\in\CS$ with an LMI:
\begin{align*}
\hat{\CK}_{\CS,\mathrm{rlx}3} &=  \{ (\SE^\top\SE)^{-1}\SE^\top (\tlZ\tlQ^{-1}) \SE
    :
    \exists\tlZ\in\BD_{\QGactive},\,
    \tlQ\in \BD_{\QGactive}^{++},\,\\
    &\rho\in\BR
    \text{ s.t. } \tlQ\tlA^\top + \tlA\tlQ + \tlZ^\top\tlB^\top + \tlB\tlZ + \rho M
     \prec O\}.
\end{align*}
This LMI can be obtained just by replacing the nonlinear term $\rho \tlQ M\tlQ$ in $\hat{\CK}_\CS$ (Theorem \ref{theorem:controller_LMI}) with the linear term $\rho M$.
Since the constraint $\tlQ M + M\tlQ \succeq \eta M,\,\eta>0$ has been removed, there is no guarantee that any $K\in \hat{\CK}_{\CS,\mathrm{rlx}3}$ is stabilizing.
Nevertheless, this approach usually achieves not only stabilization but also high numerical performance, as shown in Section \ref{sec:examples}.
\end{rem}

\subsection{$H_\infty$ control}\label{subsec:H_2_infty}

This subsection presents the $H_\infty$ control version of the results in Subsection \ref{subsec:solution_to_problem_1}.
Note that the following discussion can be applied to $H_2$ control analogously.

Consider the following perturbed system:
 \begin{align}\label{eq:LTI_system_disturbed}
         &
        \dot{x}(t) = 
        Ax(t) + Bu(t) 
         +  B_w w(t)\\
         \label{eq:LTI_system_output}
         &y(t) = Cx(t)  + Du(t) + D_w w(t),
    \end{align}
where $w$ and $y$ are an exogenous disturbance and the output, respectively.
Under $u=0$,
strict upper bounds of $H_\infty$ norm of the transfer function matrix $C(sI-A)^{-1}B_w+D_w$ are characterized by the following \textit{bounded real lemma}.
\begin{lemma}[Bounded real lemma \cite{scherer2000linear,dullerud2013course}]\label{lemma:bounded_real}
Let $G(s)=C(sI-A)^{-1}B_w+D_w$.
Then, the following statements are equivalent:
\begin{enumerate}
    \item $A$ is Hurwitz and $\|G(s)\|_\infty<\gamma$.
    \item There exists $P\succ O$ such that 
    \begin{equation}\label{eq:bounded_real_lemma}
        \begin{bmatrix}
            A^\top P+PA & PB_w & C^\top \\
            B_w^\top P   & - \gamma I & D_w^\top\\
            C & D_w & -\gamma I
        \end{bmatrix}
        \prec O.
    \end{equation}
\end{enumerate}
\end{lemma}

Consider static state feedback in \eqref{state_feedback} for the system above.
Then,
by this lemma and the notions in Section \ref{sec:preliminaries}, one can obtain an $H_\infty$ control version of Lemma \ref{lemma:stability_analysis} as follows.
For $C,\,D,\,B_w$ in \eqref{eq:LTI_system_disturbed} and \eqref{eq:LTI_system_output}, we define $\tlC$, $\tlD$, and $\tlB_w$ as
\begin{align}\label{def:extended_matrices_CD}
 \!\!\!\!\!\!
 \tlC = C(\SE^\top\SE)^{-1} \SE^\top,\,
  \tlD = D(\SE^\top\SE)^{-1} \SE^\top,\, \tlB_w = \SE B_w.
\end{align}
For an interpretation of them, see Subsection \ref{subsec:inclusion_principle}. 
\begin{lemma}\label{lemma:stability_analysis_Hinfty}
Suppose Assumption \ref{assumption:Q}.
Let 
\begin{align*}
\CP_\CS^{\infty,\gamma} =& \{P\in\CS:
    P\succ O\text{ s.t. \eqref{eq:bounded_real_lemma} holds.}
    \}, \\
\hat{\CP}_{\CS}^{\infty,\gamma} =& \{ P=\SE^\top \tlP \SE:
    \exists\tlP \in\BD_{\QGactive}^{++}
    ,\, \rho\in\BR\\
    &
    \text{s.t. } 
    \begin{bmatrix}
            \tlA^\top \tlP+\tlP\tlA + \rho M
            & \tlP \tlB_w & \Tilde{C}^\top \\
            \tlB_w^\top \tlP   & - \gamma I & D_w^\top\\
            \Tilde{C} & D_w & -\gamma I
    \end{bmatrix} \prec O
    \}.
\end{align*}
Then, $\hat{\CP}_{\CS}^{\infty,\gamma} \subset\CP_\CS^{\infty,\gamma}$ holds.
Moreover,
    $\hat{\CP}_{\CS}^{\infty,\gamma} = \CP_\CS^{\infty,\gamma}$
holds if $\CG$ is chordal and $\QGactive=\maxQG$.
\end{lemma}
\begin{proof}
Let
$\Theta = \bdiag(M, O_{n\times n}, O_{n\times n})$,    
for which we have $\Theta^\perp = \bdiag(E, I_n, I_n)$.
Then, this lemma can be shown in a similar way to Lemma \ref{lemma:sparsity} from the following relationship for \eqref{eq:bounded_real_lemma} with $P=\SE^\top\tlP\SE$:
\begin{align*}
O&\succ
\footnotesize{\begin{bmatrix}
            A^\top(\SE^\top \tlP\SE)+(\SE^\top \tlP\SE)A
            & (\SE^\top\tlP\SE) B_w & C^\top\\
            B_w^\top (\SE^\top \tlP\SE)& - \gamma I & D_w^\top\\
            C & D_w & -\gamma I
\end{bmatrix}}\\
&=
\Theta^{\perp\top}
\footnotesize{\begin{bmatrix}
            \tlA^\top \tlP+\tlP\tlA
            & \tlP \tlB_w & \Tilde{C}^\top \\
            \tlB_w^\top \tlP   & - \gamma I & D_w^\top\\
            \Tilde{C} & D_w & -\gamma I
    \end{bmatrix}}\Theta^{\perp}.
\end{align*}
Therefore, by Finsler's lemma and $\Theta^\top\Theta=\bdiag(M,O_{n\times n},O_{n\times n})$, we obtain the LMI with respect to $\tlP$, $\rho$, and $\gamma$ in $\hat{\CP}_{\CS}^{\infty,\gamma}$.
\end{proof}

The $H_\infty$ version of Theorem \ref{theorem:controller_LMI} can also be obtained as follows.
Here, let $\CK_\CS^{\infty,\gamma}$ represent all distributed control gains $K\in\CS$ with $P\in\CS$ that achieve $\|G(s)\|_\infty<\gamma$.
\begin{theorem}\label{theorem:controller_LMI_Hinfty}
Suppose Assumption \ref{assumption:Q}.
Let
\begin{align*}
&\hat{\CK}_\CS^{\infty,\gamma}= \{(\SE^\top\SE)^{-1}\SE^\top (\tlZ\tlQ^{-1}) \SE:
    \exists \tlZ\in\BD_{\QGactive},\,
    \tlQ \in\BD_{\QGactive}^{++},\\
    &\rho\in\BR
    \text{ s.t. } 
 \Gamma_\gamma(\tlQ,\tlZ) + \bdiag(\rho \tlQ M \tlQ, O_{n\times n}, O_{n\times n})
 \prec O
    \},
\end{align*}
where
\begin{align}
\label{eq:Gamma}
\!\!\!\!
    \Gamma_\gamma(\tlQ,\tlZ)= 
    \small \begin{bmatrix}
\He (\tlA\tlQ+\tlB\tlZ) 
 &  \tlB_w & \tlQ\Tilde{C}^\top + \tlZ^\top \Tilde{D}^\top  \\
    \tlB_w^\top    & - \gamma I & D_w^\top\\
    \Tilde{C}\tlQ + \Tilde{D}\tlZ & D_w & -\gamma I
    \end{bmatrix}.
\end{align}
Then, $\hat{\CK}_\CS^{\infty,\gamma} \subset \CK_\CS^{\infty,\gamma}$ holds.
Moreover,
    $\hat{\CK}_\CS^{\infty,\gamma} = \CK_\CS^{\infty,\gamma}$
holds if $\CG$ is chordal and $\QGactive=\maxQG$.
\end{theorem}
\begin{proof}
    This theorem can be shown in the same way as the proof of Theorem \ref{theorem:controller_LMI} by pre- and post-multiplying the LMI in $\hat{\CP}_\CS^{\infty,\gamma}$ with $A+BK,\,K=(\SE^\top\SE)^{-1}\SE^\top\tlK\SE$ as $A$
    by $\bdiag(\tlQ,I,I)$ with $\tlQ=\tlP^{-1}$ and setting $\tlZ=\tlK\tlQ$.
\end{proof}

A similar convex relaxation to Theorem \ref{theorem:LMI_1} can also be obtained as follows.
This theorem guarantees that the proposed relaxation in $\hat{\CK}_{\CS,\mathrm{rlx}1}^{\infty,\gamma}$ strictly covers the conventional block-diagonal relaxation.
\begin{theorem}[Proposed method 1 for $H_\infty$ control]\label{theorem:LMI_Hinfty_1}
Suppose Assumption \ref{assumption:Q}. Let 
\begin{align*}
\CK_{\CS,\diag}^{\infty,\gamma} &= \{K=ZQ^{-1}:\\
&\exists Q=\bdiag(Q_1,\ldots,Q_N)\succ O,\,Z\in\CS \text{ s.t. } \\    
&       \begin{bmatrix}
        \He (AQ+BZ) & B_w & QC^\top+Z^\top D^\top \\
            B_w^\top    & - \gamma I & D_w^\top\\
            CQ+DZ & D_w & -\gamma I
        \end{bmatrix} \prec O
    \},
\end{align*}
\begin{align*}
\hat{\CK}_{\CS,\mathrm{rlx}1}^{\infty,\gamma}=& \{(\SE^\top\SE)^{-1}\SE^\top (\tlZ\tlQ^{-1}) \SE:
    \exists \tlZ \in\BD_{\QGactive},\,
    \tlQ\in\BD_{\QGactive}^{++},\\ 
    &\rho\in\BR,\,\eta>0 \text{ s.t. }  \tlQ M + M \tlQ \succeq \eta M,\\
    & \Gamma_\gamma(\tlQ,\tlZ) + \bdiag(\rho M, O_{n\times n}, O_{n\times n})\prec O \}
\end{align*}
with $\Gamma_\gamma(\tlQ,\tlZ)$ in \eqref{eq:Gamma}.
Then, the following inclusion holds:
\begin{equation*}
   \CK_{\CS,\diag}^{\infty,\gamma} \subset\hat{\CK}_{\CS,\mathrm{rlx}1}^{\infty,\gamma}
   \subset\hat{\CK}_{\CS}^{\infty,\gamma} \subset 
   \CK_\CS^{\infty,\gamma}.
\end{equation*}
Moreover, $\CK_{\CS,\diag}^{\infty,\gamma} \subset\hat{\CK}_{\CS,\mathrm{rlx}1}^{\infty,\gamma}
   =\hat{\CK}_{\CS}^{\infty,\gamma} =
   \CK_\CS^{\infty,\gamma}$ if $\CG$ is complete.
\end{theorem}

Similarly, we can present the following variant based on the same approach as Corollary \ref{coro:LMI_2}.
\begin{coro}[Proposed method 2 for $H_\infty$ control]\label{coro:LMI_Hinfty_1}
Suppose Assumption \ref{assumption:Q}. Let 
\begin{align*}
\hat{\CK}_{\CS,\mathrm{rlx}2}^{\infty,\gamma}&= \{ K=(\SE^\top\SE)^{-1}\SE^\top (\tlZ\tlQ^{-1}) \SE:\\
    &\exists \tlZ\in\BD_{\QGactive},\,
    \tlQ \in\BD_{\QGactive}^{++}
    \text{ s.t. } 
 \Gamma_\gamma(\tlQ,\tlZ)
 \prec O
    \}
\end{align*}
with $\Gamma_\gamma(\tlQ,\tlZ)$ in \eqref{eq:Gamma}.
Then, $\hat{\CK}_{\CS,\mathrm{rlx}2}^{\infty,\gamma} \subset 
\hat{\CK}_{\CS}^{\infty,\gamma}$ holds.
Moreover, $\CK_{\CS,\diag}^{\infty,\gamma} \subset\hat{\CK}_{\CS,\mathrm{rlx}2}^{\infty,\gamma}
   =\hat{\CK}_{\CS}^{\infty,\gamma} =
   \CK_\CS^{\infty,\gamma}$ if $\CG$ is complete.
\end{coro}
\begin{rem}[Proposed method 3 for $H_\infty$ control]\label{rem:LMI_Hinfty_without_stab}
In the same spirit as $\hat{\CK}_{\CS,\mathrm{rlx}3}$ in Remark \ref{rem:LMI_without_stab}, we present the following set with an LMI, which shows interestingly strong numerical performance despite the lack of theoretical guarantee of stabilization:
\begin{align*}
&\hat{\CK}_{\CS,\mathrm{rlx}3}^{\infty,\gamma}= \{(\SE^\top\SE)^{-1}\SE^\top (\tlZ\tlQ^{-1}) \SE:
    \exists \tlZ \in\BD_{\QGactive},\,
    \tlQ\in\BD_{\QGactive}^{++},\,
    \\
    &\rho\in\BR \text{ s.t. } 
 \Gamma_\gamma(\tlQ,\tlZ) + \bdiag(\rho M, O_{n\times n}, O_{n\times n})
 \prec O\}.
\end{align*}
\end{rem}

\subsection{An interpretation of the proposed methods}\label{subsec:inclusion_principle}

We here view the proposed methods and the matrices $\tlA$, $\tlB$, $\tlC$, $\tlD$, $\tlB_w$ in \eqref{def:extended_matrices} and \eqref{def:extended_matrices_CD}, which are naturally obtained by combining Lemma \ref{lemma:agler} with Finsler's lemma, from the perspective of a state transformation with $\SE$.

Consider a LTI system in \eqref{eq:LTI_system_disturbed} and \eqref{eq:LTI_system_output} with $u=Kx$ in \eqref{state_feedback}.
For this system, transforming the state $x$ as $\Tilde{x} = \SE x$ with $\Tilde{x}(0)=\SE x(0)$ 
and plugging in $K=(\SE^\top\SE)^{-1}\SE^\top \tlK \SE$ with a $\tlK$
provide the following expanded system:
\begin{align*}
\dot{\Tilde{x}} &= \SE \dot{x} = \SE (A+BK)x + \SE B_w w \\
&= (\SE A (\SE^\top\SE)\SE^\top+\SE B(\SE^\top\SE)^{-1}\SE^\top\tlK)\SE x + \SE B_w w\\
&= (\tlA+\tlB\tlK)\Tilde{x} + \tlB_w w \\
y &= (C+DK)x + D_w w \\
&= (C(\SE^\top \SE)^{-1}\SE^\top +D(\SE^\top\SE)^{-1}\SE^\top\tlK)
\SE x + D_w w\\
&= (\tlC + \tlD \tlK) \Tilde{x} + D_w w.
\end{align*}
From Proposition \ref{prop:CD_matrix}a, one can recover $x$ from $\Tilde{x}$ by $x = (\SE^\top\SE)^{-1}\SE^\top \Tilde{x}$.

The expanded system above with $\Tilde{x} = \SE x$ tells us that
our main results and proposed methods provide a distributed controller as $\tlK\Tilde{x}$ with $\tlK=\bdiag(\ldots,\tlK_k,\ldots)$ for the expanded system of $\Tilde{x}$, and the terms $\rho M$ in $\hat{\CP}_\CS$ (Lemma \ref{lemma:stability_analysis}) and
$\rho \tlQ M \tlQ$ in $\hat{\CK}_\CS$ (Theorem \ref{theorem:controller_LMI}) are needed due to the constraint $\Tilde{x} \in \mathrm{Im}(\SE)$.
Interestingly, this state transformation was studied under the framework of \textit{inclusion principle} \cite{siljak2011decentralized,ikeda1982inclusion} in the 1980s in the control community while Agler's theorem (Lemma \ref{lemma:agler}) was initially exploited in the operations research community \cite{grone1984positive,vandenberghe2015chordal,zheng2021chordal}.

\section{Numerical examples}\label{sec:examples}

Here, we conduct numerical experiments for the stabilization and $H_\infty$ control problems
using the \texttt{mosek} \cite{andersen2000mosek}, a standard commercial SDP solver.
For the details of their implementation, see our codes available on GitHub.\footnote{
\scriptsize{\url{https://github.com/WatanabeYuto/Distributed_Control_with_Non-Block-Diagonal_Lyapunov}
}
}

\paragraph{Stabilization problem}
In the stabilization case, we set the following system and parameters.
We set $N=32$, $n_1=m_1=\ldots=n_{32}=m_{32}=1$, and 
$B =\diag(b_1,\ldots,b_{32})$ with $b_i=0,\,i=1,16$ and $b_i=1,\,i \neq 1,16$.
As $A$ matrix, we generated 200 unstable matrices such that
    $(A,B)$ are stabilizable.
    Each entry is generated by the normal Gaussian distribution.
As undirected graph $\CG$, we use 
two different graphs: ring and wheel graphs \cite{rosen2007discrete}.
\footnote{In the wheel graph, we assign node 1 to the central node.}

The simulation results are shown in Table \ref{table:stab_result}.
In the simulations, we design a stabilizing control gain for the system \eqref{eq:LTI_system} with the input in \eqref{state_feedback} by using our proposed methods 1--3, the block-diagonal relaxation, SI-based relaxation \cite{furieri2020sparsity}, and extended LMI approach \cite{ebihara2004new,pipeleers2009extended}.
In the extended LMI approach,
we set the additional scalar parameter as $\alpha=1$, similarly to \cite{ebihara2004new,ferrante2019design};
see Proposition 4 in Appendix \ref{Appendix:combination} for the details.
Table \ref{table:stab_result} indicates the number of cases where a stabilizing gain is generated out of 200 samples.

From Table \ref{table:stab_result}, our proposed methods outperform the others in all the graphs, including the extended LMI approach that generates a dense Lyapunov function.
Moreover,
the proposed method 3 always provides a stabilizing control gain despite the absence of its guarantee.
The SI-based relaxation ends in the same performance as the block-diagonal relaxation, which implies the conservativeness\footnote{
This is mainly due to the symmetry of Lyapunov matrices.
Note that \cite{furieri2020sparsity} has multiple notable benefits and insights, especially in its frequency domain formulation.}.

\renewcommand\arraystretch{1.2}
\begin{table}[t]
\centering
\scriptsize
 \caption{
 The comparison of how many cases a stabilizing control gain is obtained out of 200 different examples.
 We compare P1--P3 (proposed methods 1--3) with BD (the block-diagonal relaxation), SI (SI-based relaxation \cite{furieri2020sparsity}), and Extended LMI \cite{ebihara2004new,pipeleers2009extended}.
 }
 \label{table:stab_result}
\begin{tabular}{c||c|c|c|c|c|c}
   Graph
 &
  P1&
  P2&
  P3
  &
 BD
  &
   SI
  & 
Extended LMI
\\ \hline
Ring
& 130
& 114
& 200
& 51
& 51
& 73
\\ 
\hline
Wheel
& 149
& 178
& 200
& 54
& 54
& 106
\end{tabular}
\end{table}

\paragraph{$H_\infty$ controller design}
We also present 
numerical results of $H_\infty$ controller design
for several realistic models from $COMPle_ib$ \cite{leibfritz2004compleib},
a standard benchmark library for controller evaluation.
We remark that Appendix \ref{appendix:random_Hinf_sim} presents additional simulation results for 
randomly generated unstable $A$ matrices.
Here, we test our three proposed methods for the three models DIS1, BDT1, and DIS3 from \cite{leibfritz2004compleib}
in comparison with
the block-diagonal relaxation \eqref{LMI_diag}, SI-based relaxation \cite{furieri2020sparsity}, and the centralized $H_\infty$
controller.
The systems' dimensions are given by $(n,m)=(8,4)$, $(6,4)$, and $(11,3)$, respectively.
As the $H_\infty$ performance measure, we define
$C=[20I_n,O_{m\times m}]^\top$, $D=[O_{n\times n},200I_m]^\top$, and $D_w=O$.
As the sparsity pattern $\CS$ in the case of DIS1, 
we consider
\begin{equation*}
    P =
\tiny{\left[
\begin{array}{rrrrrrrr}
* & * & * & * & * & * & * & * \\
* & * & * & 0 & 0 & 0 & 0 & * \\
* & * & * & * & 0 & 0 & 0 & 0 \\
* & 0 & * & * & * & 0 & 0 & 0 \\
* & 0 & 0 & * & * & * & 0 & 0 \\
* & 0 & 0 & 0 & * & * & * & 0 \\
* & 0 & 0 & 0 & 0 & * & * & * \\
* & * & 0 & 0 & 0 & 0 & * & *
\end{array}
\right]}.
\end{equation*}
For the sparsity pattern of gain $K$, we take the first $4$ rows of $P$'s above.
Then, we can apply our present formulation 
by replacing the nonsquare $B$ matrix in $\mathbb{R}^{8\times 4}$ by $[B,O_{8\times 4}]\in\mathbb{R}^{8\times 8}$ without loss of generality\footnote{In fact, it is not difficult to extend our formulation to non-square matrices. For the details, see \cite[Appendix A]{fushimi2024distributed}.}.
Note that for the other examples, 
we define the sparsity pattern of $P$ and $K$ using the wheel graph in the same way.

{\renewcommand\arraystretch{1.2}
\begin{table}[t]
\centering
\scriptsize
 \caption{
 The comparison of 
 the optimal
 $H_\infty$ norm bound $\gamma^*$ for three systems (DIS1, BDT1, and DIS3) from \cite{leibfritz2004compleib}.
 The second column shows
 the optimal value $\gamma_\mathrm{cen}^*$ of the centralized controller,
 and the third to seventh columns show the value of the ratio $\gamma^*/\gamma_\mathrm{cen}^*$.
 The cases where a stabilizing controller could not be generated are denoted by "fail".}
 \label{table:Hinf_complib}
\begin{tabular}{c|c||c|c|c|c|c}
    & 
Centralized 
 &
  P1&
  P2&
  P3
  &
 BD
  &
   SI
  \\ 
\hline
DIS1
& 289.41
& 9.6728 
& 78.874 
& 1.210
& fail
& fail
\\
\hline
BDT1
& 55.702
& 1.0101
& 1.6027
& 1.0015
& 1.0504 
& 1.0504 
\\  
\hline
DIS3
& 204.886
& 1.1015
& 1.5313
& fail
& 1.1016
& 1.1016
\end{tabular}
\end{table}
}

The simulation results are presented in Table \ref{table:Hinf_complib}.
For three systems (DIS1, BDT1, and DIS3) from \cite{leibfritz2004compleib}, the second column represents
 the optimal value $\gamma_\mathrm{cen}^*$ of the centralized controller.
The third to seventh columns show the value of the ratio $\gamma^*/\gamma_\mathrm{cen}^*$, where $\gamma^*$ is the optimal $H_\infty$ norm bound by each method\footnote{
Only in the proposed method 3,
we compute $\gamma^*= \|G\|_{\infty}$ with the optimal gain $K^*$ using the MATLAB function "\texttt{hinfnorm}", due to the lack of theoretical guarantee for recovering the original bounded real lemma.}.
When a stabilizing controller could not be obtained, we denote "fail".
From this table, our proposed methods outperform the other relaxation methods, as expected from our theoretical results.
In particular, all the proposed methods achieve minimization for DIS1 while the previous methods fail.
Note that despite the failure of the proposed method 3 for DIS3, we see the highest performance in the other examples, showing its value as a heuristic algorithm.

\begin{figure}
    \centering    \includegraphics[width=0.9\linewidth]{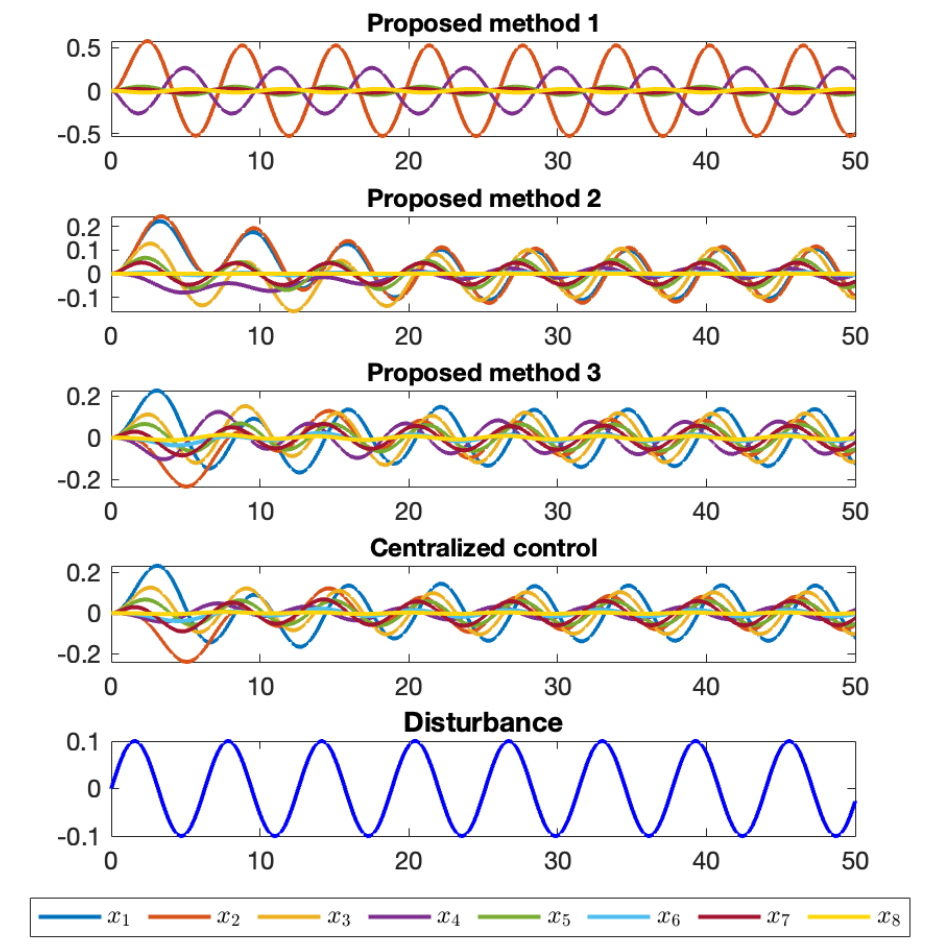}
    \caption{
    The evolution of the state vector $x(t)=[x_1(t),\ldots,x_8(t)]^\top$ in the system DIS1 \cite{leibfritz2004compleib} 
    with our proposed methods
    and the centralized $H_\infty$ controller against $w(t)=0.1\sin(t)$.}
    \label{fig:DIS1_plot}
\end{figure}
Furthermore, 
to test their disturbance attenuation performance, we also present the plots of the evolution of the states $x(t)=[x_1(t),\ldots,x_8(t)]^\top$ against $w(t) = 0.1 \sin(t)$ for DIS1 in
Fig.\ \ref{fig:DIS1_plot}.
We note that
$B_w=[1,0,1,0,1,0,1,0]^\top$ in DIS1, and we set 
the initial state as
$x(0)=0$.
This figure reveals that all the methods successfully remain bounded against the disturbance.
In particular,
the magnitude of variation due to the disturbance is also maintained small,
showing the efficacy of our proposed methods.

\section{Conclusion}\label{sec:conclusion}
This study considered the LMI-based distributed controller design problem with non-block-diagonal Lyapunov functions with the same sparsity pattern as controllers.
Based on a block-diagonal factorization of sparse matrices and Finsler's lemma,
we derived a new matrix inequality condition for distributed controllers with such Lyapunov functions and then showed that this inequality is reduced to a necessary and sufficient condition under chordal sparsity.
We further derived a new LMI as the convex relaxation and also provided analogous results for $H_\infty$ control.
Finally, numerical results demonstrated the effectiveness of our proposed methods.

\input{main.bbl}

\onecolumn
\appendix

\subsection{Proof of Lemma
\ref{lemma:sparsity}}\label{Appendix:proof-sparsity}

First, we show $\CS \supset \{\SE^\top \tlK \SE : \tlK = \bdiag(\cdots,\tlK_k,\cdots)\}$.
    The $(i,j)$ block of $\SE^\top \tlK \SE = \sum_{k \in\QGactive} \SE_{\CC_k}^\top \tlK_k \SE_{\CC_k}$
    can be represented by
    \begin{align*}
        \SE_i (\SE^\top \tlK \SE) \SE_j^\top 
        &=\sum_{k\in\QGactive} \underset{p\in\CC_k}{[\ldots,\SE_i\SE_p^\top,\ldots]} \tlK_k 
        \begin{bmatrix}
            \vdots\\
            \SE_p\SE_j^\top\\
            \vdots
        \end{bmatrix}\\
        &=\begin{cases}
            \sum_{k\in\QiGactive\cap\QjGactive} [\tlK_k]_{ij}, & \QiGactive\cap\QjGactive\neq \emptyset\\
            O, & \text{otherwise.}
        \end{cases} 
    \end{align*}
Thus, by the assumption of $\QiGactive\cap\QjGactive\neq \emptyset \Leftrightarrow (i,j)\in\CE$, we obtain $\SE^\top \tlK\SE\in\CS$ from \eqref{def:S}.
Next, we show the converse relationship.
For $K\in \CS$, consider $    \tlK_k =  \SE_{\CC_k} \hat{K}    \SE_{\CC_k}^\top$,
where $\hat{K}$'s $(i,j)$ block is defined by
\begin{equation*}
    \hat{K}_{ij} = 
    \begin{cases}
            \frac{1}{|\QiGactive\cap\QjGactive|} K_{ij}, & \QiGactive\cap\QjGactive\neq \emptyset\\
            O, & \text{otherwise.}
        \end{cases} 
\end{equation*}
Then, it holds that $\SE^\top 
    \bdiag(\cdots,\tlK_k,\cdots)
    \SE =K$.

\subsection{Proof of Proposition \ref{prop:perp}}\label{Appendix:perp}
First, we have $\SE^\top M=\SE^\top - \SE^\top=O$.
Further, for $N=\SE(\SE^\top\SE)^{-1}\SE^\top$, the definition of $\SE$ and Lemma \ref{prop:CD_matrix}a guarantee that $\mathrm{rank}(\SE)=n$, 
$\mathrm{Im}(\SE) = \mathrm{Im}(N)$, and $\mathrm{Im}(M)=\mathrm{Ker}(\SE^\top)$, where $\mathrm{Ker}(\SE^\top)=\{z:\SE^\top z = 0\}$. (For the properties of $M$ as a projection matrix, see \cite{yanai2011projection}.)
Hence, we get $\mathrm{rank}(M)= \sum_{k\in\QGactive} n_{\CC_k} - n$ and $\mathrm{rank}[M,\SE] = \sum_{k\in\QGactive} n_{\CC_k}$.
Therefore, $\SE$ can be assigned as $M^\perp$.
Notice that the assumption of $\QiGactive\neq\emptyset$ ensures $\sum_{k\in\QGactive} n_{\CC_k}\geq n$.
\subsection{A supporting result for Theorem \ref{theorem:LMI_1}}

Here, we present a supporting result for proving Theorem \ref{theorem:LMI_1}.
This proposition is key to showing
the superiority of our proposed method 1 to the block-diagonal relaxation.
\begin{prop}\label{prop:diag-Q-E}
    Consider the same assumptions as Theorem \ref{theorem:LMI_1}.
    Let $Q = \bdiag(Q_1,\ldots,Q_N)\succ O$ with $Q_i\in\BR^{n_i\times n_i}$.
    Then, for $\tlQ =\bdiag(\cdots,\tlQ_k,\cdots)\succ O$ with
        $\tlQ_k = \bdiag(\cdots,\underset{j\in\CC_k}{
        |\QjGactive| Q_j}, \cdots)\succ O$, 
    we have
    \begin{equation*}
        M\tlQ = \tlQ M.
    \end{equation*}
\end{prop}
\begin{proof}
Since the definition of $\tlQ$ gives $Q^{-1} = E^\top \tlQ^{-1}E$,
we obtain
\begin{align*}
\tlQ M 
&= \tlQ - \tlQ E(E^\top E)^{-1}E^\top \\
&= \tlQ - 
\begin{bmatrix}
\vdots\\
\tlQ_k E_{\CC_k}\\
\vdots
\end{bmatrix}(E^\top E)^{-1}E^\top \\
&= \tlQ - 
\underbrace{\begin{bmatrix}
\vdots\\
E_{\CC_k}(E^\top E)Q\\
\vdots
\end{bmatrix}}_{=E(E^\top E)Q=EQ(E^\top E)}(E^\top E)^{-1}E^\top \\
&=\tlQ - EQE^\top.
\end{align*}
As $\tlQ - EQE^\top$ is symmetric, $\tlQ M = (\tlQ M)^\top = M\tlQ$.
\end{proof}

\subsection{Combining the proposed approach with early works}\label{Appendix:combination}

Here, 
we demonstrate how to combine our proposed approach with existing studies and how the performance is improved.
As a result, we establish the inclusion relationships in Fig. \ref{fig:appendix_incliusion_diagram} with new LMIs, where we define $\CK$ as
 \begin{equation}\label{def:K_all}
     \CK :=
     \{K\in\CS:
     \exists P\succ O
     \text{ s.t. }
\eqref{eq:original_problem}
     \}.
 \end{equation}
We also provide an additional numerical result; see Remark \ref{rem:num_combined} and Table \ref{table:stab_result_aux}.
 
The following results are built upon our trick below using $P=\SE^\top \tlP \SE$ and $K=(\SE^\top\SE)^{-1}\SE^\top \tlK\SE$, which is also a key to our main results in Lemma \ref{lemma:stability_analysis} and Theorem \ref{theorem:controller_LMI}:
\begin{align*}
&(A+BK)^\top P + P(A+BK)\\
&= \SE^\top ((\tlA+\tlB\tlK)^\top \tlP + \tlP(\tlA+\tlB\tlK))\SE \prec O\\
\overset{\text{Lem \ref{lemma:finsler}}}{\Leftrightarrow}
&\;
\exists\rho\in\BR\text{ s.t. }
(\tlA+\tlB\tlK)^\top \tlP + \tlP(\tlA+\tlB\tlK) + \rho M \prec O\\
\Leftrightarrow&
\;
\exists\rho\in\BR\text{ s.t. }
\tlA\tlQ+\tlQ\tlA^\top + \tlB\tlZ+ \tlZ^\top\tlB
+\rho \tlQ M\tlQ\prec O.
\end{align*}

Note that for the LMI in \cite{yuan2023structured}, we have not shown an explicit inclusion relationship as Fig. \ref{fig:appendix_incliusion_diagram}, but we show that our approach is helpful to identify the source of conservatism and allows one to search control gains that cannot be obtained by \cite{yuan2022parameter,yuan2023structured}.

\begin{figure}[t]
\centering
\begin{tikzpicture}
\definecolor{setcolor1}{RGB}{254, 178, 76}
\definecolor{setcolor2}{RGB}{199, 233, 180}
\definecolor{setcolor3}{RGB}{127, 205, 187}
\definecolor{setcolor4}{RGB}{65, 182, 196}

\draw[fill=setcolor1, fill opacity=0.1] (0,0) ellipse (4 and 3.5);
\node at (0, 3.2) {\(\CK\)\text{ in \eqref{def:K_all}}};

\draw[fill=setcolor4, fill opacity=0.5] (-0.5, 0) ellipse (3.2 and 3);
\node at (-2.5, 2) {\(
\substack{\hat{\CK}_{\CS,\mathrm{rlx}1}^{\mathrm{FDR}} \\
\text{ in \eqref{def:K_FDR}}}
\)
};
\draw[fill=setcolor4, fill opacity=0.5] (0.5, 0) ellipse (3.2 and 3);
\node at (2.5, 2) {\(
\substack{\hat{\CK}_{\CS,\mathrm{rlx}1}^{\mathrm{FZPK}} \\
\text{ in \eqref{def:K_FZPK}}}
\)
};
\draw[fill=white, fill opacity=0.8] (0,0) ellipse (2.5 and 2.5);
\node at (0, 2) {\(
\substack{
\hat{\CK}_{\CS,\mathrm{rlx}1}^{\mathrm{ext}}\\
\text{in Theorem \ref{theorem:K_ext_rlx}}
}
\)};

\draw[fill=setcolor1, fill opacity=0.3] (-0.5,-0.2) ellipse (1.8 and 1.8);
\node at (-1.3, 1.0) {\(
\substack{{\CK}_{\CS,\mathrm{ext}}\\
\text{in \eqref{def:K_ext}}
}
\)};

\draw[fill=setcolor1, fill opacity=0.3] (0.5,-0.2) ellipse (1.8 and 1.8);
\node at (1.3, 1.0) {\(\substack{
{\hat{\CK}}_{\CS,\mathrm{rlx}1}\\
\text{in Theorem \ref{theorem:LMI_1}}
}
\)};

\draw[fill=white, fill opacity=1] (0,-0.2 ) ellipse (1 and 1);
\node at (0, 0) {\(\CK_{\CS,\mathrm{diag}}\text{\scriptsize{ in \eqref{LMI_diag}}}\)};

\end{tikzpicture}
\caption{
The inclusion relationships between the block-diagonal relaxation $\CK_{\CS,\mathrm{diag}}$,
proposed method 1 $\hat{\CK}_{\CS,\mathrm{rlx}1}$ in Theorem \ref{theorem:LMI_1},
extended LMI ${\CK}_{\CS,\mathrm{ext}}$ in \eqref{def:K_ext},
combined method
$\hat{\CK}_{\CS,\mathrm{rlx}1}^{\mathrm{ext}}$ with the proposed method and extended LMI in Theorem \ref{theorem:K_ext_rlx}, and the other combinations
$\hat{\CK}_{\CS,\mathrm{rlx}1}^{\mathrm{FDR}}$
and 
$\hat{\CK}_{\CS,\mathrm{rlx}1}^{\mathrm{FZPK}}$.
}
\label{fig:appendix_incliusion_diagram}
\end{figure}

{
\renewcommand\arraystretch{1.2}
\begin{table*}
[t]\label{table:additional-result_stability}
\small
\centering
 \caption{
 The comparison of how many cases a stabilizing control gain is obtained out of 200 different examples in Table. \ref{table:stab_result} with the combined method.
 We compare P1--P3 (proposed methods 1--3) 
 and the combined method in Theorem \ref{theorem:K_ext_rlx}
 with BD (the block-diagonal relaxation), SI (SI-based relaxation \cite{furieri2020sparsity}), and Extended LMI \cite{ebihara2004new,pipeleers2009extended}.
 In the combined method and extended LMI, we set $\alpha=1$.
 }
 \label{table:stab_result_aux}
\begin{tabular}{c||c|c|c|c|c|c|c}
   Graph
 &
  P1&
  P2&
  P3
  &
 BD
  &
   SI
  & 
  \begin{tabular}{c}
       Extended\\
       LMI
  \end{tabular}&
  Combined method
  \\ \hline
Ring
& 130
& 114
& 200
& 51
& 51
& 73
& 200
\\ 
\hline
Wheel
& 149
& 178
& 200
& 54
& 54
& 106
& 199
\end{tabular}
\end{table*}
}

\paragraph{Extended LMI techniques}
Here, we show that our proposed approach $\hat{\CK}_{\CS,\mathrm{rlx}1}$ can be combined with the well-known extended LMI approach (denoted by $\CK_{\CS,\mathrm{ext}}$ Fig. \ref{fig:appendix_incliusion_diagram}), which can improve standard LMIs by introducing a slack variable.
The resultant LMI (denoted by in $\hat{\CK}^{\mathrm{ext}}_{\CS,\mathrm{rlx}1}$ Fig. \ref{fig:appendix_incliusion_diagram})
provides a broader class of $K\in\CK_\CS$ than the two methods.

Note that this extended LMI technique using a slack variable was first established for discrete-time systems, and
not all but some results have been extended to continuous-time systems.
The following formulations are the continuous-time versions of the authors' recent paper \cite{fushimi2024distributed} dedicated to discrete-time systems.

As a preliminary, we illustrate the basic idea of the extended LMI technique.
In the stabilizing controller design via extended LMI,
the following proposition is fundamental.
Here, let $H$ represent the slack variable.
\begin{prop}[\cite{pipeleers2009extended,ebihara2004new,ferrante2019design}]\label{prop:extended_LMI}
For the LTI system in \eqref{eq:LTI_system} with a state-feedback controller $u=Kx$ in \eqref{state_feedback} and an arbitrarily prescribed number $\alpha>0$,
the following items are equivalent:
\begin{enumerate}
    \item $A+BK$ is stable.
    \item $P\succ  O$ and $P(A+BK)+(A+BK)^\top P \prec O$ hold. 
    \item $P\succ  O$ and
\begin{equation}\label{eq:extended_LMI}
     \begin{bmatrix}
         O & P\\
         P & O
     \end{bmatrix}
     + \He\left(
      \begin{bmatrix}
          (A+BK)^\top\\
          -I
      \end{bmatrix}
      H
      \begin{bmatrix}
          I & \alpha I
      \end{bmatrix}
     \right)
     \prec O
    \end{equation}
    with a square matrix $H$ hold.
\end{enumerate}
\end{prop}

Equivalently transforming \eqref{eq:extended_LMI} into
\begin{equation}\label{eq:extended_LMI_transformed}
     \begin{bmatrix}
         O & Q\\
         Q & O
     \end{bmatrix}
     + \He\left(
      \begin{bmatrix}
          G^\top A^\top +Z^\top B^\top \\
          -G^\top
      \end{bmatrix}
      \begin{bmatrix}
          I & \alpha I
      \end{bmatrix}
     \right)
     \prec O
    \end{equation}
with the new variables as $G=H^{-1}$, $Q:=(H^{-1})^\top PH^{-1}$,
$Z = KG$ with $Z \in \CS$, and $G=\bdiag(G_1,\ldots,G_N)$
provides a convex reformulation of \eqref{def:LMI_nonconv},
by which we obtain a stabilizing distributed controller $K=ZG^{-1}\in\CS$.
Remarkably this approach enables dispensing with the sparsity constraint on $P$ matrix since the constraint is imposed on $G$ instead.
Namely,
we can obtain a distributed control gain as a member of the following set:
\begin{align}\label{def:K_ext}
\CK_{\CS,\mathrm{ext}}
=\{&K = ZG^{-1}:
\exists Q\succ O,\,Z\in\CS,
\nonumber
\\
&G=\bdiag(G_1,\ldots,G_N) \text{ s.t. } \eqref{eq:extended_LMI_transformed}
\}.
\end{align}
It is easy to see that
\begin{equation}\label{eq:diag-ext}
\CK_{\CS,\mathrm{diag}} \subset \CK_{\CS,\mathrm{ext}} \subset \CK_\CS
\end{equation}
because for $P=H$, \eqref{eq:extended_LMI} is reduced to $P(A+BK)+(A+BK)^\top P\prec O.$
Note that the non-singularity of $G$ is guaranteed as $G+G^\top \succ 0$ from \eqref{eq:extended_LMI_transformed} (see the (2,2) block).

Now, we show the following theorem, which demonstrates
that combining our proposed approach with the extended LMI presents a novel set $\hat{\CK}_{\CS}^{\mathrm{ext}}$.
Notice that
setting $\Tilde{G}=\Tilde{Q}$ in this set yields $\hat{\CK}_{\CS}$, which implies
$\hat{\CK}_{\CS}\subset\hat{\CK}_{\CS}^{\mathrm{ext}}$.

\begin{theorem}\label{theorem:K_ext}
Suppose Assumption \ref{assumption:Q}.
Let $\alpha>0$ and
\begin{align*}
\hat{\CK}_{\CS}^{\mathrm{ext}} 
=
\{&K=(\SE^\top\SE)^{-1}\SE^\top \tlZ\Tilde{G}^{-1}\SE:\exists\tlQ\succ O,\,
 \tlZ = \bdiag(\tlZ_1,\ldots,\tlZ_q),\,\Tilde{G} = \bdiag(\Tilde{G}_1,\ldots,\Tilde{G}_q),\, \rho\in\BR  \\
&
 \text{ s.t. }\begin{bmatrix}
         O & \tlQ\\
         \tlQ & O
 \end{bmatrix}
     + \He\left(
      \begin{bmatrix}
          \Tilde{G}^\top \tlA^\top +\tlZ^\top \tlB^\top \\
          -\Tilde{G}^\top
      \end{bmatrix}
      \begin{bmatrix}
          I & \alpha I
      \end{bmatrix}
     \right)+ \bdiag(\rho\Tilde{G} M\Tilde{G},\rho\Tilde{G} M\Tilde{G})
     \prec O
\}.
\end{align*}
Then, $\hat{\CK}_\CS \subset\hat{\CK}_{\CS}^{\mathrm{ext}}\subset \CK$ holds.
\end{theorem}
\begin{proof}
    One can show this theorem similarly to the discrete-time version in \cite{fushimi2024distributed}.
    The proof sketch is as follows.
    First, we have $\hat{\CK}_\CS \subset\hat{\CK}_{\CS}^{\mathrm{ext}}$ as we recover the matrix inequality in $\hat{\CK}_\CS$ 
    from that in $\hat{\CK}_{\CS}^{\mathrm{ext}}$
    by $\Tilde{G}=\tlQ$.
    The latter inclusion $\hat{\CK}_{\CS}^{\mathrm{ext}}\subset\CK$ can also be verified by pre- and post-multiplying the matrix inequality in $\hat{\CK}_{\CS}^{\mathrm{ext}}$ by $E^\top \Tilde{G}^{-1}$ and $\Tilde{G}^{-1}E$ and
    applying Finsler's lemma.
\end{proof}

We can apply
the same relaxation technique as Theorem \ref{theorem:LMI_1} to
the nonlinear matrix inequality in $\hat{\CK}_{\CS}^{\mathrm{ext}}$,
which gives a novel LMI
in the set $\hat{\CK}_{\CS,\mathrm{rlx}1}^{\mathrm{ext}}$ below.
As shown in Fig. \ref{fig:appendix_incliusion_diagram},
this set covers both the extended LMI approach $\CK_{\CS,\mathrm{ext}}$ and proposed method 1 $\hat{\CK}_{\CS,\mathrm{rlx}1}^{\mathrm{ext}}$.

\begin{theorem}[Combined method with the extended LMI]\label{theorem:K_ext_rlx}
Suppose Assumption \ref{assumption:Q}.
Let $\alpha>0$ and
\begin{align*}
\hat{\CK}_{\CS,\mathrm{rlx}1}^{\mathrm{ext}} 
=
\{&K=(\SE^\top\SE)^{-1}\SE^\top \tlZ\Tilde{G}^{-1}\SE:\exists\tlQ\succ O,\,
 \tlZ = \bdiag(\tlZ_1,\ldots,\tlZ_q),\,\Tilde{G} = \bdiag(\Tilde{G}_1,\ldots,\Tilde{G}_q),\, \rho\in\BR,\,\eta>0 \\
&\text{ s.t. } 
 \begin{bmatrix}
         O & \tlQ\\
         \tlQ & O
 \end{bmatrix}
     + \He\left(
      \begin{bmatrix}
          \Tilde{G}^\top \tlA^\top +\tlZ^\top \tlB^\top \\
          -\Tilde{G}^\top
      \end{bmatrix}
      \begin{bmatrix}
          I & \alpha I
      \end{bmatrix}
     \right)+ \bdiag(\rho M,\rho M)
     \prec O,\, \Tilde{G}^\top M + M\Tilde{G} \succeq \eta M
\}.
\end{align*}
Then, 
\begin{align*}
&
\CK_{\CS,\mathrm{diag}}\subset
\CK_{\CS,\mathrm{ext}} 
\subset \hat{\CK}_{\CS,\mathrm{rlx}1}^{\mathrm{ext}}
\subset
\hat{\CK}_{\CS}^{\mathrm{ext}}
\subset\CK \\
&\CK_{\CS,\mathrm{diag}}\subset
\hat{\CK}_{\CS,\mathrm{rlx}1} 
\subset 
\hat{\CK}_{\CS,\mathrm{rlx}1}^{\mathrm{ext}}\subset
\hat{\CK}_{\CS}^{\mathrm{ext}}
\subset
\CK.
\end{align*}
\end{theorem}
\begin{proof}
We first show $\CK_{\CS,\mathrm{ext}} 
\subset \hat{\CK}_{\CS,\mathrm{rlx}1}^{\mathrm{ext}}$ as follows.
Consider $K\in\CK_{\CS,\mathrm{ext}}$.
Then there exist $Z\in\CS$ and $G=\bdiag(G_1,\ldots,G_N)$ such that $K= ZG^{-1}$.
For $G$, let us define 
$\Tilde{G} = \bdiag(\Tilde{G}_1,\ldots,\Tilde{G}_q)$
with 
    $\Tilde{G}_k = \bdiag(\cdots,\underset{j\in\CC_k}{|\QjGactive|G_j},\cdots)$.
For this matrix, $\Tilde{G}M=M\Tilde{G}$ 
and $E^\top \Tilde{G}^{-1} E  =G^{-1}$
hold from the same reasoning as Proposition \ref{prop:diag-Q-E}.
Then, invoking $G+G^\top \succ 0$, we have
\begin{align*}
\Tilde{G}^\top M+M\Tilde{G}&=\Tilde{G}^\top M^2+M^2\Tilde{G}= M(\Tilde{G}^\top +\Tilde{G})M \succeq \eta M
\end{align*}
with $\eta = \lambda_\mathrm{min}(\Tilde{G}^\top +\Tilde{G})>0$.
Moreover, since there exists $\tlK\in\BD_{\QGactive}$ satisfying $K=(E^\top E)^{-1}E^\top \tlK E$ from Lemma \ref{lemma:sparsity}, it can be seen that
\begin{align*}
    O \prec  &
    \begin{bmatrix}
        (G^{\top})^{-1}& (G^{\top})^{-1}
    \end{bmatrix}
    \biggl(
    \begin{bmatrix}
         O & Q\\
         Q & O
     \end{bmatrix}+ \He\left(
      \begin{bmatrix}
          G^\top A^\top +Z^\top B^\top \\
          -G^\top
      \end{bmatrix}
      \begin{bmatrix}
          I & \alpha I
      \end{bmatrix}
     \right)
     \biggr)
     \begin{bmatrix}
         G^{-1}\\G^{-1}
     \end{bmatrix}
     \\
     =& 
        \begin{bmatrix}
         O & (G^{\top})^{-1}QG^{-1}\\
         (G^{\top})^{-1}QG^{-1} & O
     \end{bmatrix}+ \He\left(
      \begin{bmatrix}
          (A+BK)^\top\\
          -I
      \end{bmatrix}
      G^{-1}
      \begin{bmatrix}
          I & \alpha I
      \end{bmatrix}
     \right) \\
     =& 
        \begin{bmatrix}
         O & (G^{\top})^{-1}QG^{-1}\\
         (G^{\top})^{-1}QG^{-1} & O
     \end{bmatrix}+ \He\left(
      \begin{bmatrix}
          E^\top (\tlA+\tlB \tlK)^\top  \\
          -E^\top 
      \end{bmatrix}
      \Tilde{G}^{-1}
      \begin{bmatrix}
          E & \alpha E
      \end{bmatrix}
     \right).
     \end{align*}
 This is further equivalent to 
     \begin{align*}
     \Leftrightarrow \quad&
     \exists \nu \in\BR \text{ s.t. } \\
    &\begin{bmatrix}
         O & \tlP \\
         \tlP & O
     \end{bmatrix}    
     + \He\left(
      \begin{bmatrix}
          (\tlA+\tlB \tlK)^\top  \\
          -I
      \end{bmatrix}
      \Tilde{G}^{-1}
      \begin{bmatrix}
          I & \alpha I
      \end{bmatrix}
     \right) + \bdiag(\nu M ,\nu M)\prec O\\
      \Leftrightarrow \quad&
     \exists \nu \in\BR \text{ s.t. } \\
    &\begin{bmatrix}
         O & \tlQ \\
         \tlQ& O
     \end{bmatrix}    
     + \He\left(
      \begin{bmatrix}
          \Tilde{G}^\top \tlA^\top +\tlZ^\top \tlB^\top \\
          - \Tilde{G}^\top
      \end{bmatrix}
      \begin{bmatrix}
          I & \alpha I
      \end{bmatrix}
     \right) + \bdiag(\nu  \Tilde{G}^\top M \Tilde{G} ,\nu \Tilde{G}^\top M \Tilde{G})\prec O,
\end{align*}
where $\tlP$ is a positive definite matrix satisfying $E^\top \tlP E = (G^{\top})^{-1}QG^{-1}$, and we define $\tlQ=\Tilde{G}^\top \tlP\Tilde{G}\succ O$
and 
$\tlZ = \tlK \Tilde{G} \in \BD_{\QGactive}$.
Since it holds for the last term that
\begin{equation*}
    -\nu \Tilde{G}^\top M \Tilde{G}
    =
    -\nu \Tilde{G}^\top  M^2 \Tilde{G}
    \preceq |\nu| M\tilde{G}^\top\Tilde{G} M
    \preceq -\rho M
\end{equation*}
with $\rho =- |\nu| \lambda_\mathrm{max}(\Tilde{G}^\top\Tilde{G})$,
we arrive at 
\begin{align*}
    &\begin{bmatrix}
         O & \tlQ \\
         \tlQ& O
     \end{bmatrix}    
     + \He\left(
      \begin{bmatrix}
          \Tilde{G}^\top \tlA^\top +\tlZ^\top \tlB^\top \\
          - \Tilde{G}^\top
      \end{bmatrix}
      \begin{bmatrix}
          I & \alpha I
      \end{bmatrix}
     \right) + \bdiag(\rho M ,\rho M)\prec O, \,\Tilde{G}^\top M+M\Tilde{G} \succeq \eta M
\end{align*}
with $\tlQ\succ 0$, $\tlZ \in \BD_{\QGactive}$, $\Tilde{G}\in\BD_{\QGactive}^{++}$,
$\rho\in\BR$, and $\eta>0$,
which implies $K \in \hat{\CK}_{\CS,\mathrm{rlx}1}^\mathrm{ext}$.
Therefore, $\mathcal{K}_{\CS,\mathrm{ext}}\subset \hat{\CK}_{\CS,\mathrm{rlx}1}^\mathrm{ext}$.

Finally,
the other parts can be obtained in a similar manner to our main results;
the inclusion $\hat{\CK}_{\CS,\mathrm{rlx}1}^{\mathrm{ext}}\subset
\hat{\CK}_{\CS}^{\mathrm{ext}}$ follows from the same argument as the proof of Theorem \ref{theorem:LMI_1} based on Finsler's lemma.
It is also straightforward to show $\hat{\CK}_{\CS,\mathrm{rlx}1} 
\subset \hat{\CK}_{\CS,\mathrm{rlx}1}^{\mathrm{ext}}$ because
the set $\hat{\CK}_{\CS,\mathrm{rlx}1} $ is obtained by setting $\Tilde{G}=\Tilde{Q}$ in $\hat{\CK}_{\CS,\mathrm{rlx}1}^{\mathrm{ext}}$.
The rest of the inclusions are the immediate consequences of Theorems \ref{theorem:controller_LMI}, \ref{theorem:LMI_1}, and \ref{theorem:K_ext} and the inclusion \eqref{eq:diag-ext}.
\end{proof}

\begin{rem}\label{rem:num_combined}
    In Table \ref{table:stab_result_aux},
    we present additional numerical results compared with \cite{pipeleers2009extended} and the combined method in Theorem \ref{theorem:K_ext_rlx}
    for two different graphs.
    Similar to the numerical simulation in Section \ref{sec:examples},
    we set $\alpha=1$ for the combined method and extended LMI.
    It can be observed that the combined method achieves highly competitive performance with the guarantee of yielding a stabilizing gain.
\end{rem}

\paragraph{Ferrante et al. (2019) \cite{ferrante2019design} and Furieri et al. (2020) \cite{furieri2020sparsity}}

Leveraging the decomposition techniques for the constraint $Z Q^{-1}\in\CS$ in Ferrante et al. (2019) \cite{ferrante2019design} and in Furieri et al. (2020) \cite{furieri2020sparsity},
we can obtain the following two sets with LMIs.
Here, we apply the same techniques to $\tlZ \tlQ^{-1}\in\bdiag(\tlK_1,\ldots,\tlK_q)$ in the proposed approach instead of $Z Q^{-1}\in\CS$.
Note that the following subscripts of FDR and FZPK represent the authors' initials.
\begin{align}\label{def:K_FDR}
\begin{aligned}
\hat{\CK}_{\CS,\mathrm{rlx}1}^{\mathrm{FDR}} 
=
\{&K=(\SE^\top\SE)^{-1}\SE^\top \tlZ\Tilde{G}^{-1}\SE:\exists\tlQ\succ O,\,
 \tlZ \in \mathrm{span}\{S_1,\ldots,S_{r}\},\,
\Tilde{G} \in \Upsilon
,\, \rho\in\BR,\,\eta>0\\
&
 \text{ s.t. } 
 \begin{bmatrix}
         O & \tlQ\\
         \tlQ & O
 \end{bmatrix}
     + \He\left(
      \begin{bmatrix}
          \Tilde{G}^\top \tlA^\top +\tlZ^\top \tlB^\top \\
          -\Tilde{G}^\top
      \end{bmatrix}
      \begin{bmatrix}
          I & \alpha I
      \end{bmatrix}
     \right)+ \bdiag(\rho M,\rho M)
     \prec O,\, \Tilde{G}^\top M + M\Tilde{G} \succeq \eta M
\},
\end{aligned}
\end{align}
\begin{align}
\begin{aligned}\label{def:K_FZPK}
\hat{\CK}_{\CS,\mathrm{rlx}1}^{\mathrm{FZPK}} 
=
\{&K=(\SE^\top\SE)^{-1}\SE^\top \tlZ\Tilde{G}^{-1}\SE:\exists\tlQ\succ O,\,
 \tlZ \in \Tilde{\CS}_Z,\,
\Tilde{G}\in \Tilde{\CS}_G,
,\, \rho\in\BR,\,\eta>0  \\
&
\text{ s.t. }
 \begin{bmatrix}
         O & \tlQ\\
         \tlQ & O
 \end{bmatrix}
     + \He\left(
      \begin{bmatrix}
          \Tilde{G}^\top \tlA^\top +\tlZ^\top \tlB^\top \\
          -\Tilde{G}^\top
      \end{bmatrix}
      \begin{bmatrix}
          I & \alpha I
      \end{bmatrix}
     \right)+ \bdiag(\rho M,\rho M)
     \prec O,\, \Tilde{G}^\top M + M\Tilde{G} \succeq \eta M
\}\end{aligned}
\end{align}
where $\{S_1,\ldots,S_r\}$ represent the basis of $\CS$,
\begin{align*}
\Upsilon
:=
\{V: \exists \Lambda \in \mathcal{R}^r
\text{ s.t. }
L(I\otimes V) = L(\Lambda \otimes I)
\}
\end{align*}
with $L := [S_1|S_2|\ldots|S_r]$, and the set of non-singular matrices $\mathcal{R}^r$.
Further, the sets $\Tilde{\CS}_Z$ and $\Tilde{\CS}_G$ satisfy
\begin{equation*}
\tlZ\in\Tilde{\CS}_Z,\,
\Tilde{G}\in\Tilde{\CS}_G
\Rightarrow
\tlZ\Tilde{G}^{-1}\in\Tilde{\CS},
\end{equation*}
where $\Tilde{\CS}= \{\tlK: \tlK= \bdiag(\tlK_1,\ldots,\tlK_q)\}$.
For their definitions and computation details, see \cite{ferrante2019design} and \cite{furieri2020sparsity}.

As illustrated in Fig. \ref{fig:appendix_incliusion_diagram}, we have the following results for these sets:
\begin{align*}
&
\hat{\CK}_{\CS,\mathrm{rlx}1}^{\mathrm{ext}}
\subset
\hat{\CK}_{\CS,\mathrm{rlx}1}^{\mathrm{FDR}}
\subset\CK \\
&\hat{\CK}_{\CS,\mathrm{rlx}1}^{\mathrm{ext}}
\subset
\hat{\CK}_{\CS,\mathrm{rlx}1}^{\mathrm{FZPK}}
\subset\CK.
\end{align*}
These inclusions can easily be verified by setting $\tlZ = \bdiag(\tlZ_1,\ldots,\tlZ_q)$ and
$\Tilde{G} = \bdiag(\Tilde{G}_1,\ldots,\Tilde{G}_q)$ in $\hat{\CK}_{\CS,\mathrm{rlx}1}^{\mathrm{FDR}}$ and $\hat{\CK}_{\CS,\mathrm{rlx}1}^{\mathrm{FZPK}}$.



\paragraph{Yuan et al. (2023) \cite{yuan2023structured}}
These works
present an LMI for distributed $H_2$ controller design using a similar block diagonal factorization.
Here, we show that 
this approach is built upon a similar idea to our proposed method 2 in Corollary \ref{coro:LMI_2}, and
our proposed method 1 provides a different solution set of stabilizing control gains.

In \cite{yuan2023structured},
the following convex restriction of the $H_2$ optimal controller design problem 
has been presented based on \cite{yuan2022parameter}:
\begin{align}\label{problem:Yuan2023}
\begin{array}{cl}
\underset{
\substack{\tlQ=\bdiag(\tlQ_1,\ldots,\tlQ_q),\\
\tlZ=\bdiag(\tlZ_1,\ldots,\tlZ_q),\,\Tilde{W}
}
}
{\mbox{minimize}}&
\displaystyle  \mathrm{trace}(\tlQ\tlC^\top\tlC + \Tilde{W} D_w^\top D_w)
\\
{\mbox{subject to}}&  \displaystyle 
\He(\tlA\tlQ+\tlB\tlZ) + \tlB_w\tlB_w^\top \preceq O,\\
& \tlQ\succ O,\\
& \begin{bmatrix}
    \tlQ & \tlZ^\top \\
    \tlZ & \Tilde{W}
\end{bmatrix} \succeq O.
\end{array}
\end{align}
Solving this problem yields a distributed controller $K=(\SE^\top\SE)^{-1}\SE^\top\tlK\SE$ with $\tlK=\tlZ\tlQ^{-1}$.
Here, the closed-loop stability is guaranteed by the first and second constraints as follows:
\begin{align}\label{eq:stability_yuan}
    \begin{aligned}
            &\SE^\top \tlQ^{-1}(\He(\tlA\tlQ+\tlB\tlZ) + \tlB_w\tlB_w^\top)\tlQ^{-1}\SE \\
    =&\SE^\top (\He(\tlP\tlA+\tlP\tlB\tlK) + \tlP\tlB_w\tlB_w^\top\tlP)\SE\\
    =& \He(P(A+BK)) + PB_wB_w^\top P \preceq O,
    \end{aligned}
\end{align}
where $\tlP=\tlQ^{-1}$ and
$P=\SE^\top\tlP\SE$.

Notice that Problem \eqref{problem:Yuan2023} causes conservatism because the transformation in \eqref{eq:stability_yuan} still holds when we replace the first constraint in \eqref{problem:Yuan2023} by
the following nonlinear matrix inequality with some $\rho\in\BR$ as Theorems \ref{theorem:controller_LMI} and \ref{theorem:controller_LMI_Hinfty}:
\begin{equation}\label{eq:QMQ_yuan}
\He(\tlA\tlQ+\tlB\tlZ) + \tlB_w\tlB_w^\top 
+\rho \tlQ M\tlQ
\preceq O.
\end{equation}
Namely, our result reveals that
Problem \eqref{problem:Yuan2023} in \cite{yuan2023structured} corresponds to the case of $\rho = 0$ in \eqref{eq:QMQ_yuan} as the proposed method 2, which can be conservative because
if \eqref{eq:QMQ_yuan} holds for $\rho=0$, so does for all $\rho <0$ from $\tlQ M \tlQ\succeq O$.
In other words,  Problem \eqref{problem:Yuan2023} in \cite{yuan2023structured} can be quite conservative when the optimal value of $\rho$ is negative.
Note that if $\CG$ is chordal,
\eqref{eq:QMQ_yuan} provides a lossless expression, albeit nonlinearity, which is another insight from our approach.

To mitigate the conservatism of Problem \eqref{problem:Yuan2023} and overcome the nonlinearity of \eqref{eq:QMQ_yuan},
one can apply the same relaxation technique as Theorems \ref{theorem:LMI_1} and \ref{theorem:LMI_Hinfty_1}, which allows us to explore the cases where the optimal value of $\rho$ is negative.

\begin{rem}
    Note that
    in the proposed method 2, we use the strict inequality
    to ensure the closed-loop stability
    while the authors of \cite{yuan2023structured} employ the non-strict inequality in the $H_2$ control problem as \eqref{problem:Yuan2023}, which does not lead to instability thanks to the additional term, $\tlB_w\tlB_w^\top$.
\end{rem}
\begin{rem}\label{remark:yuan}
    We remark that 
    although our proposed approach sheds new light on a missing perspective in \cite{yuan2023structured}, this method cannot be derived from our proposed approach based on Finsler's lemma in the same way, since it 
    directly applies 
    the relaxation 
$\tilde{Q}=\bdiag(\tilde{Q}_1,\ldots,\tilde{Q}_q)$ and $\tilde{K}=
\bdiag(\tilde{X}_1,\ldots,\tilde{X}_q)
$ 
to a \textit{dual} LMI in the higher-dimensional space of $\mathrm{Im}(E)$, not in $\mathbb{R}^n$, 
    and thus does not emerge from
    primal LMIs of the form
    $\He(P(A+BK))\prec O$.
    Note also that in \cite{yuan2023structured}, the authors present a highly flexible block-diagonal factorization beyond Lemmas \ref{lemma:agler} and \ref{lemma:sparsity}.
    Our approach can also be combined with this factorization.
\end{rem}

\subsection{Extensions of the proposed approach to systems with polytopic uncertainties and static output feedback}\label{Appendix:extention}

Here, we discuss the way to deal with polytopic uncertainties and static output feedback
through the proposed approach.
Thanks to the systematic nature of our approach,
these extensions are not very difficult as follows.

\paragraph{Polytopic uncertainties}

Polytopic linear systems can be treated through our approach in the same way
because 
each of the matrices in \eqref{def:extended_matrices} is obtained by a linear transformation for the system matrices of the original system in \eqref{eq:LTI_system_disturbed} and \eqref{eq:LTI_system_output}.
For example,
when we have
\begin{equation*}
    A = \sum_{p=1}^P \alpha_p A_p
\end{equation*} with $\alpha_p\in[0,1],\,p=1,\ldots,P$ satisfying
$\sum_{p=1}^P\alpha_p = 1$,
then it holds that
\begin{align*}
     \tlA &= \SE(\sum_{p=1}^P \alpha_p A_p)(\SE^\top\SE)^{-1}\SE^\top =
     \sum_{p=1}^P \alpha_p 
     (\SE A_p(\SE^\top\SE)^{-1}\SE^\top ).
\end{align*}
From the perspective of the state transformation with $\SE$ in Subsection \ref{subsec:inclusion_principle},
polytopic uncertainty of the original system of $x$ results in that of the expanded system of $\tilde{x}$.
Hence, thanks to this property, 
our approach is still effective even in the face of polytopic uncertainty.

\paragraph{Static output-feedback problem}
Our approach can also be extended to static output-feedback.
Consider the following system:
\begin{align*}
    &\dot{x} = Ax + Bu\\
    &y = Cx
\end{align*}
with a row full rank matrix $C$.
For designing a static output feedback controller, the following well-known LMI was proposed in the late 1990s in \cite{crusius1999sufficient} as a sufficient condition:
\begin{equation}\label{eq:LMI-output-fb}
\begin{array}{l}AQ+QA^{\top}+B Z C+C^{\top} Z^\top B^{\top} \prec O \\ Q\succ O \\ W C=C Q,\end{array}
\end{equation}
which provides the output feedback controller $u = Ky$ with
\begin{equation*}
    K = Z W^{-1}
\end{equation*}
since $C Q^{-1} = W^{-1}C$ holds.
A typical distributed output feedback controller design using \eqref{eq:LMI-output-fb} is as follows.
Consider that the output $y$ is given by
\begin{equation*}
  y = \begin{bmatrix}
      y_1\\
      \vdots\\
      v_N
  \end{bmatrix}
  = 
  \begin{bmatrix}
      C_1 x_1\\
      \vdots\\
      C_N x_n
  \end{bmatrix}
  =   \bdiag(C_1,\ldots,C_N)x.
\end{equation*}
Then, a distributed controller $u = Ky$ with $K=ZW^{-1} \in \CS$ is obtained
by solving the following LMI for $Q$, $Z$, and $W$:
\begin{align}
\label{LMI:output-fb-distributed}
    \begin{array}{l}AQ+QA^{\top}+B Z C+C^{\top} Z^\top B^{\top} \prec O \\ Q =\bdiag(Q_1,\ldots,Q_N) \succ O \\ W C=C Q\\
    Z\in\CS.
    \end{array}
\end{align}
This LMI can be viewed as \eqref{eq:LMI-output-fb} with the block-diagonal relaxation.
Note that we have to adjust the size of $K$ and $Z$ when $C$ is not square.

Our proposed approach allows us to
improve the LMI above.
Suppose that Assumption \ref{assumption:K_sparsity} is satisfied and matrix $C$ is row full rank.
Then,
applying the same approach as Theorem \ref{theorem:LMI_1} to \eqref{eq:LMI-output-fb} provides
the following LMI and a distributed output feedback controller $u = (\SE^\top\SE)^{-1}\SE^\top \tlK y$ with $\tlK=\tlZ W^{-1}$:
\begin{equation*}
      \begin{array}{l}
      \tlA\tlQ+\tlQ\tlA^{\top}+\tlB \tlZ \tlC + \tlC^{\top} \tlZ^\top \tlB^{\top} 
      + \rho M \prec O\\
      \tlQ =\bdiag(\tlQ_1,\ldots,\tlQ_q) \succ O \\ W \tlC=\tlC \tlQ\\
    \tlZ = \bdiag(\tlZ_1,\ldots,\tlZ_q) 
    \prec O \\
  \tlQ M  + M \tlQ \succeq \eta M\\
  \rho \in \BR,\, \eta >0
    \end{array}
\end{equation*}
as $\tlC = C(\SE^\top\SE)^{-1}\SE^\top$ is also column full rank under Assumption \ref{assumption:K_sparsity}.
Notice that the original LMI in \eqref{LMI:output-fb-distributed}
can be recovered
in the same way
as the proof of Theorem \ref{theorem:LMI_1}.

\begin{figure*}
\begin{subfigure}{.5\textwidth}
  \centering
\includegraphics[width=0.85\columnwidth]{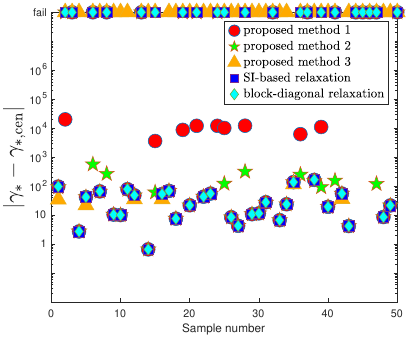}
    \caption{The simulation results for the ring graph.}
    \label{fig:results_ring}
\end{subfigure}%
\begin{subfigure}{.5\textwidth}
  \centering
\includegraphics[width=0.85\columnwidth]{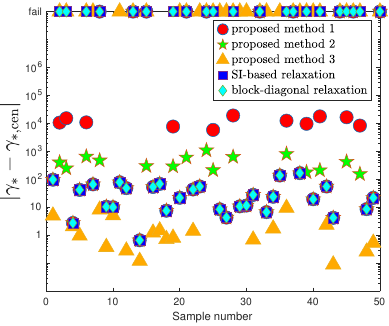}
    \caption{The simulation results for the wheel graph.}
    \label{fig:results_disk}
\end{subfigure}
\caption{The red circles, stars, orange circles, squares, and diamonds represent the value of $|\gamma_*-\gamma_{*,\mathrm{cen}}|$ under the proposed method 1 in $\CK_{\CS,\mathrm{rlx}1}^{\infty,\gamma}$, proposed method 2 in $\hat{\CK}_{\CS,\mathrm{rlx}2}^{\infty,\gamma}$,
    proposed method 3 in $\hat{\CK}_{\CS,\mathrm{rlx}3}^{\infty,\gamma}$,
    SI-based relaxation \cite{furieri2020sparsity}, and block-diagonal relaxation in $\CK_{\CS,\diag}^{\infty,\gamma}$, respectively.
    Here, $\gamma_{*,\mathrm{cen}}$ represents the optimal value of the upper bound $\gamma$ of $\|(C+DK)(sI-(A+BK))^{-1}B_w + D_w\|_\infty$ under the centralized controller.
    The cases where a stabilizing controller was not generated are collected at the top as failures.}
\label{fig:hinf_sim}
\end{figure*}

\subsection{Additional numerical results for $H_\infty$ control}\label{appendix:random_Hinf_sim}

This section presents 
additional simulation results for distributed $H_\infty$ controller design using randomly generated systems.
We used the solver \texttt{mosek} in the simulation results below.

We here minimize $\gamma$, i.e., the upper bound of the $\|(C+DK)(sI-(A+BK))^{-1}B_w + D_w\|_\infty$ for the system in \eqref{eq:LTI_system_disturbed} and \eqref{eq:LTI_system_output} with a distributed controller \eqref{state_feedback} with $K\in\CS$. 
Here, we let $B=\diag(0,1,\ldots,1)$ and randomly generate 50 $A$ matrices that have 6 unstable eigenvalues and satisfy the stabilizability.
As an undirected graph, we employ the same two graphs as the stabilization case, i.e., ring and wheel graphs.
As the $H_\infty$ performance measure, we set 
$C=[20I_n,O_{m\times m}]^\top$, $D=[O_{n\times n},I_m]^\top$, and $D_w=O$.

The simulation results for $50$ different random samples for the two graphs
are plotted in Figs \ref{fig:results_ring} and \ref{fig:results_disk}, respectively.
The horizontal axis represents the sample number of each pair of random $A$ and $\CG$, and the vertical axis corresponds to $|\gamma_*-\gamma_{*,\mathrm{cen}}|$, 
where $\gamma_{*,\mathrm{cen}}$ and $\gamma_*$ represent the optimal value of $\gamma$ for the centralized controller and for the others, respectively\footnote{
Only in the proposed method 3,
we compute $\gamma= \|G\|_{\infty}$ with the optimal gain $K$, due to the lack of theoretical guarantee for recovering the original bounded real lemma.}.
The red circles, stars, orange circles, squares, and diamonds represent the value of $|\gamma_*-\gamma_{*,\mathrm{cen}}|$ under the proposed method 1 in $\CK_{\CS,\mathrm{rlx}1}^{\infty,\gamma}$,
proposed method 2 in $\CK_{\CS,\mathrm{rlx}2}^{\infty,\gamma}$, proposed method 3 in $\hat{\CK}_{\CS,\mathrm{rlx}3}^{\infty,\gamma}$, SI-based relaxation \cite{furieri2020sparsity}, and block-diagonal relaxation in $\CK_{\CS,\diag}^{\infty,\gamma}$, respectively.
The cases where the generated distributed controller $K$ did not stabilize the system due to numerical problems are collected at the top of Figs \ref{fig:results_ring} and \ref{fig:results_disk} as failures.
We do not consider \cite{ferrante2020lmi} and
the extended LMI approach for $H_\infty$ control \cite{shaked2001improved}
for fairness because they involve non-convexity (e.g., a bilinear scalar parameter) and require an additional step for convexification.
Note that for the discrete-time case, 
we provide a fair numerical and theoretical comparison in \cite{fushimi2024distributed}.

From Figs \ref{fig:results_ring} and \ref{fig:results_disk}, 
the proposed methods 1 and 2 provide a feasible solution most frequently, demonstrated by the largest number of generating of a stabilizing gain, and achieve superior or equivalent performance to the other previous methods. Further, the proposed method 3 successfully generates a stabilizing gain in many cases with the smallest $\gamma$ for the wheel graph that is chordal, albeit in the absence of its theoretical guarantee.
In both graphs (including the non-chordal ring one), the best results with a closer performance to the centralized controller can be observed in the proposed method 3 when it succeeds.

\renewcommand\arraystretch{1.2}
\begin{table}[t]
\centering
 \caption{
 The comparison of how many cases a stabilizing $H_\infty$ suboptimal gain with $\gamma=1000$ is obtained out of the same 50 different examples as Fig. \ref{fig:hinf_sim}.
 We compare P1--P3 (proposed methods 1--3) with SI (SI-based relaxation \cite{furieri2020sparsity}) and BD (the block-diagonal relaxation).
For the performance of the proposed method 3, see Remark \ref{remark:opt_vs_subopt}.
 }
 \label{table:hinf_suboptimal}
\begin{tabular}{c||c|c|c|c|c|c}
   Graph
 &
  P1&
  P2&
  P3
  &
   SI
    &
   BD
\\ \hline
Ring
& 40
& 46
& 50
& 29
& 29
\\ 
\hline
Wheel
& 31
& 49
& 50
& 29
& 29
\end{tabular}
\end{table}

Furthermore, we also present the simulation results of $H_\infty$ suboptimal control with $\gamma=1000$ for the same randomly generated systems.
Similarly, we compare our proposed methods with the SI-based and block-diagonal relaxations.
The simulation result is presented in Table \ref{table:hinf_suboptimal}, which shows the number of instances where a feasible stabilizing controller is successfully obtained.
From the table, the proposed method 3 achieved the best performance, followed by the other proposed methods.
From Tables \ref{table:stab_result} and \ref{table:hinf_suboptimal}, it can be observed that the proposed method 3 can offer preferable performance for feasibility problems.
These results demonstrate the efficacy of our proposed methods 1 and 2 and the value of proposed method 3 as a heuristic strategy that has complementary benefits to the other proposed methods.

\begin{rem}\label{remark:opt_vs_subopt}
    In the proposed method 3,
there is a considerable gap between minimizing $\gamma$ and finding a feasible controller with a fixed $\gamma$.
One of the possible contributing factors to the discrepancy is the loss of \textit{strict feasibility} in the dual problem when $\gamma$ is minimized, leading to the violation of the Slater condition,
a necessary condition for solvers based on the primal-dual interior point method.
Such an issue could be rooted in the $H_\infty$ \textit{optimal} control problem itself, as pointed out in several works, e.g., \cite{waki2015application}.
Due to the lack of theoretical guarantees, 
there is room for further investigation on the behavior of the proposed method 3, which is included in our future work.
\end{rem}

\end{document}